\documentclass[runningheads]{llncs}
\usepackage[T1]{fontenc}
\usepackage{amsmath,amssymb,mathtools}
\usepackage{graphicx}
\usepackage{xcolor}
\usepackage{tikz}
\usetikzlibrary{positioning,arrows.meta,calc, decorations.pathmorphing}
\usepackage{subcaption}
\usepackage[hidelinks]{hyperref}
\usepackage[capitalise,noabbrev]{cleveref}

\definecolor{ccolor1}{HTML}{28009F}
\definecolor{ccolor2}{HTML}{9F0028}
\tikzset{solidnode/.style={circle,fill=black,minimum size=1ex,inner sep=0pt},
  dedge/.style={-{Latex[width=1.1ex,length=1.1ex]},line width=.15ex},
  demandedge/.style={dedge,dashed}, wavy/.style={dedge,decorate,decoration={snake,amplitude=.4mm,
  segment length=2.5mm,post length=1.5mm}}}
\newcommand{\dist}{\mathrm{dist}}
\newcommand{\supp}{\mathrm{supp}}

\newcommand{\ind}[1]{\mathbf 1_{#1}}
\newcommand{\R}{\mathbb R}

\begin{document}
\title{Cut Conditions, Integrality, and Obstructions\\ for Directed Multiflows}
\titlerunning{Cut Conditions and Directed Multiflows}
\author{
Srinwanti Debgupta\inst{1}\and
Martin Skutella\inst{1}\and
Philipp Warode\inst{2}}
\authorrunning{S. Debgupta, M. Skutella, Ph. Warode
}
\institute{
Technische Universit\"at Berlin, Germany\\
\email{\{debgupta,skutella\}@math.tu-berlin.de}
\and Humboldt-Universit\"at zu Berlin, Germany\\
\email{philipp.warode@hu-berlin.de}
}
\maketitle
\setcounter{footnote}{0}

\begin{abstract}
In contrast to the max-flow/min-cut theorem for single-com\-modity flows, the cut condition is necessary but not sufficient for the existence of directed multiflows, even with only two commodities. Poremba and Shepherd (IPCO\,2023) conjecture that the cut condition is sufficient for all capacities and demands on a directed two-commodity topology if and only if it contains neither of two specific topologies as a relevant minor. We prove this conjecture. As a consequence, the fractional and integral flow-cut gaps of every directed two-commodity topology coincide and are either $1$ or~$2$.

We then investigate stronger necessary conditions for multiflow feasibility and their implications for integrality. We introduce a hierarchy of five conditions, ranging from the classical cut condition through conditions based on blocking edge sets and binary edge lengths to the linear programming dual characterization of fractional feasibility. The corresponding classes of sufficient topologies are strictly nested. We extend the Nagamochi--Ibaraki theorem: whenever one of the first four conditions is sufficient on a topology, integral capacities and demands satisfying that condition admit an integral multiflow, computable in polynomial time. Finally, we prove that recognizing topologies on which any of these four conditions is sufficient is coNP-complete, even for two commodities.

\keywords{directed multiflows \and cut-sufficiency \and integral flows \and relevant minors \and flow-cut gaps}
\end{abstract}

\section{Introduction}
\label{sec:introduction}

Multiflows model the simultaneous routing of several commodities
through a network with shared capacities, with applications in
communication networks, transportation, and logistics~\cite{Ahuja_Magnanti_Orlin}. Their
expressive power also connects them closely to general linear
programming: fractional multiflow feasibility has a poly\-nomial-size
linear programming formulation, and conversely, even directed
two-commodity flow feasibility is as hard as general linear
programming~\cite{Itai_78,Ding_Kyng_Zhang_22}.

A \emph{multiflow topology}~$(G,H)$ consists of a supply digraph~$G=(V,E)$
and a demand digraph~$H$ on the same vertex set.
Each demand edge~$(s_i,t_i)$ represents a commodity~$i$ with
source~$s_i$ and sink~$t_i$. A \emph{multiflow} assigns nonnegative
flow values to directed~$s_i$--$t_i$~paths for each commodity~$i$.
Given nonnegative edge capacities~$u_e$ on supply edges~$e\in E$ and demands~$d_i$, it is
\emph{feasible} if the paths of each commodity~$i$ carry a total
flow of~$d_i$ and the total flow of all commodities on each supply
edge~$e$ is at most~$u_e$. A multiflow is \emph{integral} if all
path-flow values are integers; unless stated otherwise, we allow
fractional flow values.

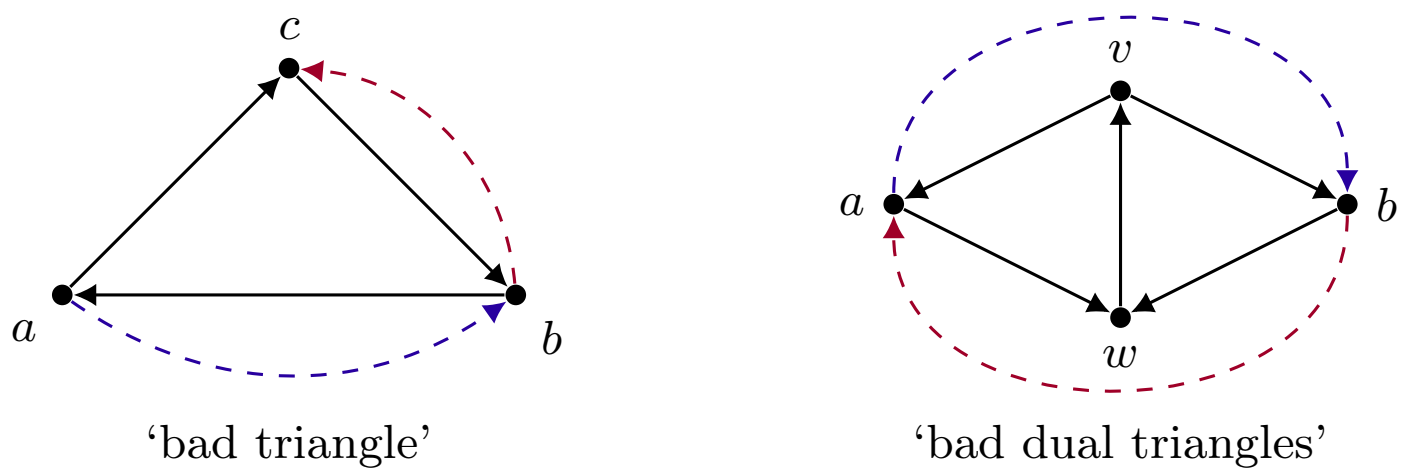
\begin{figure}[!t]
\centering
\begin{tikzpicture}[x=1cm,y=1cm]
  \begin{scope}[yshift=-0.6cm]    
    \node[solidnode,label=below left:$a$] (ta) at (0,0) {};
    \node[solidnode,label=above:$c$] (tc) at (1.5,1.5) {};
    \node[solidnode,label=below right:$b$] (tb) at (3,0) {};
    \draw (ta) edge[dedge] (tc) (tc) edge[dedge] (tb) (tb) edge[dedge] (ta);
    \draw[demandedge,ccolor1] (ta) to[bend right=35] (tb);
    \draw[demandedge,ccolor2] (tb) to[bend right=40] (tc);
  \end{scope}
  \node[font=\small] at (1.5,-1.6) {`bad triangle'};
  \begin{scope}[xshift=5.5cm]
  \node[solidnode,label=left:$a$] (da) at (0,0) {};
  \node[solidnode,label=below:$w$] (dw) at (1.5,-.75) {};
  \node[solidnode,label=above:$v$] (dv) at (1.5,.75) {};
  \node[solidnode,label=right:$b$] (db) at (3,0) {};
  \draw (da) edge[dedge] (dw) (db) edge[dedge] (dw) (dw) edge[dedge] (dv)
        (dv) edge[dedge] (da) (dv) edge[dedge] (db);
  \draw[demandedge,ccolor1] (da) .. controls (0,1.6) and (3,1.6) .. (db);
  \draw[demandedge,ccolor2] (db) .. controls (3,-1.6) and (0,-1.6) .. (da);
  \node[font=\small] at (1.5,-1.6) {`bad dual triangles'};
  \end{scope}
\end{tikzpicture}
\caption{The bad triangle (left) and the bad dual triangles (right).
Solid edges are supply edges; the dashed demand edges represent commodities. With all capacities
and demands equal to one, both topologies satisfy the cut condition but
admit no feasible multiflow: the unique paths of the two commodities
share a unit-capacity edge, namely~$(a,c)$ on the left and~$(w,v)$
on the right.}
\label{char:obstructions}
\end{figure}

For the special case of a single commodity, the max-flow/min-cut theorem characterizes
feasibility by capacity inequalities for cuts, and integral capacities
and demands guarantee an integral flow whenever a feasible flow
exists. For multiple commodities, the analogous \emph{cut condition}
requires every vertex set to have enough outgoing capacity for the
total demand of commodities whose source lies inside and sink
outside. This condition remains necessary; however, it is no longer
sufficient, even for two commodities (see
\cref{char:obstructions}). Moreover, fractional feasibility need
not imply integral feasibility for integral capacities and demands,
again already with two commodities.

The distinction between directed and undirected multiflows is also
apparent with only two commodities. In undirected graphs, the cut
condition suffices for fractional feasibility~\cite{Hu_63}.
Disjoint-path problems exhibit a similar contrast: deciding whether
a digraph contains vertex-disjoint directed paths for two given
terminal pairs is NP-hard~\cite{Fortune_Hopcroft_Wyllie_80},
whereas in undirected graphs the corresponding problem is
polynomial-time solvable for any fixed number of terminal
pairs~\cite{Robertson_Seymour_95}.

These differences motivate the study of directed topologies on which
the cut condition nevertheless characterizes feasibility. A topology
is \emph{cut-sufficient} if, for every choice of nonnegative capacities
and demands, the cut condition implies the existence of a feasible
multiflow. Nagamochi and Ibaraki~\cite{Nagamochi_Ibaraki_89} show
that cut-sufficiency also restores integrality: on a cut-sufficient
topology, integral capacities and demands satisfying the cut condition
admit an integral multiflow.

Poremba and Shepherd~\cite{Poremba_Shepherd_2025} initiated the study
of cut-sufficiency for directed two-commodity topologies in terms of
forbidden obstructions. Cut-sufficiency is not inherited by arbitrary
minors, but it is inherited by \emph{relevant minors}. A minor is
relevant if every choice of capacities and demands satisfying the
cut condition on the minor still satisfies it when extended to the
original topology: surviving edges retain their weights, deleted
supply and demand edges receive capacity and demand zero,
respectively, and contracted supply edges receive sufficiently
large (or even infinite) capacities.

Poremba and Shepherd's basic obstructions are the \emph{bad triangle} and the
\emph{bad dual triangles} (\cref{char:obstructions}).
They prove that excluding both as relevant minors
characterizes cut-sufficiency when the two demands form a directed 
path ($t_1=s_2$) or a roundtrip ($t_1=s_2$ and $t_2=s_1$); they also establish
this characterization for directed star and triangle demand graphs.
They also prove that recognizing cut-sufficiency is
coNP-hard, even for roundtrip demands. The remaining non-trivial
two-commodity case consists of four distinct terminals, called
\emph{two-matching demands}. For this case, they conjecture the same
forbidden-minor characterization and prove it for unit demands and
integral capacities~\cite[Thm.~4]{Poremba_Shepherd_2025}.

\paragraph{Main result.}
We prove the conjecture of Poremba and Shepherd (\cref{char:main}): a
directed two-commodity topology is cut-sufficient if and only if it
contains neither the bad triangle nor the bad dual triangles as a
relevant minor. In particular, sufficiency for unit demands and integral
capacities on a two-matching topology implies sufficiency for arbitrary
nonnegative capacities and demands.
The argument of Poremba and Shepherd for unit demands analyzes pairs of
paths directly. Our proof instead starts from a linear programming
certificate of infeasibility: a nonnegative integral edge-length
function that violates the dual length inequality.
From such a certificate we obtain a path from one commodity's source to
the other commodity's sink that is shorter than the first commodity's
own source--sink distance. Contracting this path yields a relevant minor
whose two demands form a directed path, the case solved by Poremba and
Shepherd. Minimality of the certificate guarantees that, after the
contraction, the two commodities still have simple paths sharing an
edge. Their characterization then supplies a forbidden minor.

As a consequence, non-cut-sufficiency of a
directed two-commodity topology is always certified by demands
in~$\{0,1\}$ and capacities in~$\{0,1,2\}$, confirming the unit-demand
certificate property conjectured by Poremba and Shepherd
(\cref{char:unit-certificate}).
Moreover, the fractional and integral flow-cut gaps of every directed
two-commodity topology coincide and equal one or two (\cref{char:gap}).
The flow-cut gap is the least factor, at least one, by which capacities
must be scaled to guarantee a routing for every choice of capacities
and demands satisfying the cut condition; in the integral case, the
original weights and the routing are required to be integral. The
upper bound two is immediate for two commodities; the content lies in
the dichotomy and in the equality of the fractional and integral gaps.
Our result determines the gap exactly for each topology. For the broader
study of fractional and integral flow-cut gaps, see Chekuri, Shepherd,
and Weibel~\cite{Chekuri_Shepherd_Weibel_13}. When the number of
commodities is unbounded, directed
flow-cut gaps can be polynomially large~\cite{Chuzhoy_Khanna_09}.

\paragraph{A hierarchy of cut conditions.}
We place these results in a hierarchy of five necessary conditions. A
commodity is \emph{blocked} by an edge set~$F\subseteq E$ if each of its paths uses an
edge in~$F$. Beyond the cut condition~\eqref{C1}, condition~\eqref{C2}
requires every cut~$F=\delta_G^+(X)$ leaving a vertex set~$X\subseteq V$
to have capacity at least the total demand it blocks.
Condition~\eqref{C3} imposes the same requirement on every edge
set~$F\subseteq E$. Both were studied by Majthoub
Almoghrabi et al.~\cite{majthoub2026integer}; \eqref{C3} is the classical
condition behind directed sparsest
cut~\cite{Harvey_Kleinberg_Lehman_06,Chuzhoy_Khanna_09}, and a framework
related to~\eqref{C2} appears in~\cite{Chekuri_Kannan_Raja_Viswanath_11}.
Condition~\eqref{C4} also counts how many edges of the set each path must
use, via $0/1$ edge lengths. Condition~\eqref{C5} allows all nonnegative
integral lengths and characterizes feasibility by linear programming
duality. The classes~$\mathcal S_r$ of topologies on which
the $r$th condition suffices are strictly nested, already for two
commodities (\cref{hier:strict}). Checking any of~\eqref{C2},
\eqref{C3}, or~\eqref{C4} is strongly coNP-complete even for two
commodities, whereas checking~\eqref{C1} for a fixed number of
commodities or~\eqref{C5} for any number takes polynomial time
(\cref{hier:easy-checks,hier:hard-checks}). Recognizing each of~$\mathcal S_1,\dots,\mathcal S_4$ is coNP-complete,
even for two-matching demands (\cref{hier:recognition}).

Extending Nagamochi and
Ibaraki's theorem~\cite{Nagamochi_Ibaraki_89}, we show that sufficiency at any of
the first four levels implies integral sufficiency, with an integral multiflow computable
in polynomial time (\cref{int:main}); for~\eqref{C2} and~\eqref{C3},
the existence statement was observed in~\cite{majthoub2026integer}. The
key property is that these conditions have $0/1$ capacity coefficients
and integral right-hand sides for integral demands.
Condition~\eqref{C4} is the strongest condition in the hierarchy with
$0/1$ capacity coefficients. The integrality conclusion holds through
this level but fails for~\eqref{C5}, where fractional feasibility need
not imply integral feasibility even for integral capacities and demands.

\paragraph{Organization.}~\Cref{sec:preliminaries} fixes notation, and~\cref{hier:sec} introduces
the hierarchy. \Cref{char:section} proves the main result and its
consequences, and~\cref{int:section} proves the integrality theorem. The appendix contains
sufficient families for the intermediate levels, including a closure
rule under vertex sums, the verification of the strict inclusions, the
complexity results, a weight-preserving contraction for blocking-set
violations, results on gaps at the higher levels, and length tests with
bounded coefficients.

\section{Preliminaries}
\label{sec:preliminaries}

All graphs are finite directed multigraphs: parallel edges are distinct, and all paths are simple directed paths. A \emph{multiflow
topology} $(G,H)$ has a supply digraph $G=(V,E)$ and a demand
digraph $H$ on the same vertex set. The $k$ demand edges (\emph{commodities}) are $(s_i,t_i)$ for
$i\in[k]:=\{1,\ldots,k\}$, identified with their indices. A
\emph{weighting} (or \emph{weights}) of~$(G,H)$ consists of
capacities~$u\in\R_{\ge0}^{E}$ and demands $d\in\R_{\ge0}^{k}$, forming
 a \emph{multiflow instance} $(G,H,u,d)$, 
\emph{integral} if $u$ and $d$ are.  Demand loops, whose source and sink
coincide, are discarded. Supply loops may also be discarded, since they
lie on no simple path and cross no cut. Zero demands and capacities are
allowed: the topology is fixed; only the weighting varies.

Let~$\mathcal P_i$ be the set of $s_i$--$t_i$ paths in~$G$. A
\emph{multiflow}~$x$ assigns a value~$x^i_P\ge0$ to every
path~$P\in\mathcal P_i$, $i\in[k]$; the total flow of commodity~$i$ on an
edge~$e$ is $x^i_e:=\sum_{P\in\mathcal P_i:\,e\in P}x^i_P$. The multiflow is
\emph{feasible} if
\[
 \sum_{P\in\mathcal P_i}x^i_P=d_i\quad\text{for }i\in[k]
 \qquad\text{and}\qquad
 \sum_{i\in[k]}x^i_e\leq u_e\quad\text{for }e\in E,
\]
and an instance is \emph{feasible} if it admits a feasible multiflow.
Multiflows are fractional unless stated otherwise; an \emph{integer
multiflow} has all~$x^i_P$ integral. Cancelling cycles shows that the usual edge-based formulation, used as a linear program in \cref{int:section}, has the same feasible instances.

For a digraph~$Q$ on~$V$ and $X\subseteq V$, let~$\delta_Q^+(X)$ and~$\delta_Q^-(X)$ denote its edges leaving and entering~$X$, respectively.
A \emph{cut} of~$G$ is an edge set~$\delta_G^+(X)$.
For $F\subseteq E$
and $I\subseteq[k]$, write $u(F):=\sum_{e\in F}u_e$ and
$d(I):=\sum_{i\in I}d_i$.

A vertex set~$X \subseteq V$ \emph{separates} commodity~$i$ if $s_i\in X$ and
$t_i\notin X$; thus~$\delta^+_H(X)$ is the set of separated commodities.
An edge set~$F\subseteq E$ \emph{blocks} commodity~$i$ if every
$s_i$--$t_i$ path uses an edge in~$F$. Let~$K_G(F)$ (or~$K(F)$ when~$G$ is clear) be
 the commodities blocked by~$F$. Every separated commodity is blocked by the corresponding cut, but a cut may also block a commodity it does not separate. To see that every blocked commodity is separated by a cut contained in~$F$, let~$X_i(F)$ be the vertices reachable from~$s_i$ in~$G-F$. Then
\begin{equation}
 s_i\in X_i(F),\qquad \delta^+_G(X_i(F))\subseteq F,
 \qquad\text{and}\qquad t_i\notin X_i(F)\iff i\in K(F).
 \label{eq:canonical-cut}
\end{equation}

For edge lengths~$\ell\in\R_{\ge0}^{E}$, write $\ell(P):=\sum_{e\in P}\ell_e$
for a path or walk~$P$, and let~$\dist_\ell(v,w)$ be the minimum length
of a $v$--$w$ path, or~$+\infty$ if none exists; in demand-weighted sums
we use $0\cdot(+\infty)=0$. Let $\supp(\ell):=\{e\in E:\ell_e>0\}$, and
let~$\ind{S}$ be the incidence vector of a set~$S$ of edges, vertices or commodities. For
vertices~$a$ before~$b$ on a path~$P$, $P[a,b]$ is the subpath from~$a$
to~$b$, and~$\cdot$ denotes concatenation of paths and edges sharing an
endpoint.
\section{A hierarchy of cut conditions}
\label{hier:sec}

The five necessary conditions are:
\begin{align}
 u(\delta_G^+(X))&\ge d(\delta_H^+(X))
 &&\text{for every }X\subseteq V,\tag{C1}\label{C1}\\
 u(\delta_G^+(X))&\ge d(K(\delta_G^+(X)))
 &&\text{for every }X\subseteq V,\tag{C2}\label{C2}\\
 u(F)&\ge d(K(F))
 &&\text{for every }F\subseteq E(G),\tag{C3}\label{C3}\\
 \sum_{e\in E(G)}u_e\ell_e
 &\ge\sum_{i=1}^k d_i\dist_\ell(s_i,t_i)
 &&\text{for every }\ell\in\{0,1\}^{E(G)},\tag{C4}\label{C4}\\
 \sum_{e\in E(G)}u_e\ell_e
 &\ge\sum_{i=1}^k d_i\dist_\ell(s_i,t_i)
 &&\text{for every }\ell\in\mathbb Z_{\ge0}^{E(G)}.\tag{C5}\label{C5}
\end{align}

Weights satisfying~\eqref{C1} in an infeasible instance are
\emph{deceptive}; a topology is cut-sufficient exactly when it has no
deceptive weights. Conditions~\eqref{C2} and~\eqref{C3} were studied
in~\cite{majthoub2026integer}, where~\eqref{C2} is shown to suffice for
series-parallel digraphs. Condition~\eqref{C3} is also the condition
underlying directed sparsest cut~\cite{Harvey_Kleinberg_Lehman_06,Chuzhoy_Khanna_09};
a framework related to~\eqref{C2} appears in
\cite[Sect.~2.2]{Chekuri_Kannan_Raja_Viswanath_11}.

For a fixed edge set~$F$, routing each commodity on a path minimizing
its number of edges in~$F$ shows that the minimum~$u(F)$ over feasible
weights with demands~$d$ is
$\sum_i d_i\dist_{\ind{F}}(s_i,t_i)$, provided every positive-demand
commodity has a path. Therefore, any inequality $u(F)\ge\sum_i q_id_i$ valid for all feasible instances
must satisfy $q_i\le\dist_{\ind{F}}(s_i,t_i)$, and is implied by~\eqref{C4} with $\ell=\ind{F}$.

Every feasible instance satisfies all five conditions.  Conversely,
\eqref{C5} implies feasibility by~\cite{Iri_71,Onaga_Kakusho_71}, or equivalently by linear
programming duality; see, e.g., \cite[Thm.~70.1]{schrijver2003combinatorial}.%
\footnote{The theorem is usually stated for real lengths. If a
positive-demand pair is unreachable, $\ell=0$ already
violates~\eqref{C5}. Otherwise, a strict violation by real lengths
persists under a sufficiently close rational approximation, and scaling
gives integral lengths.}
Moreover,
\[
 \eqref{C5}\ \Longrightarrow\ \eqref{C4}\ \Longrightarrow\
 \eqref{C3}\ \Longrightarrow\ \eqref{C2}\ \Longrightarrow\ \eqref{C1}.
\]
For~$\eqref{C4}\Rightarrow\eqref{C3}$, take $\ell=\ind{F}$.
The other implications follow from definitions.

\begin{definition}
For $r\in\{1,\ldots,5\}$, a topology~$(G,H)$ is
\emph{\textup{(C$r$)}-sufficient} if every weighting $(u,d)$
satisfying~\textup{(C$r$)} is feasible. Let~$\mathcal S_r$ be the class
of these topologies.
\end{definition}

Thus,
$\mathcal S_1\subseteq\mathcal S_2\subseteq\dots\subseteq\mathcal S_5$.
Here~$\mathcal S_1$ denotes cut-sufficient topologies,
whereas $\mathcal S_5$ contains all topologies.  All inclusions are
strict, already for two commodities.

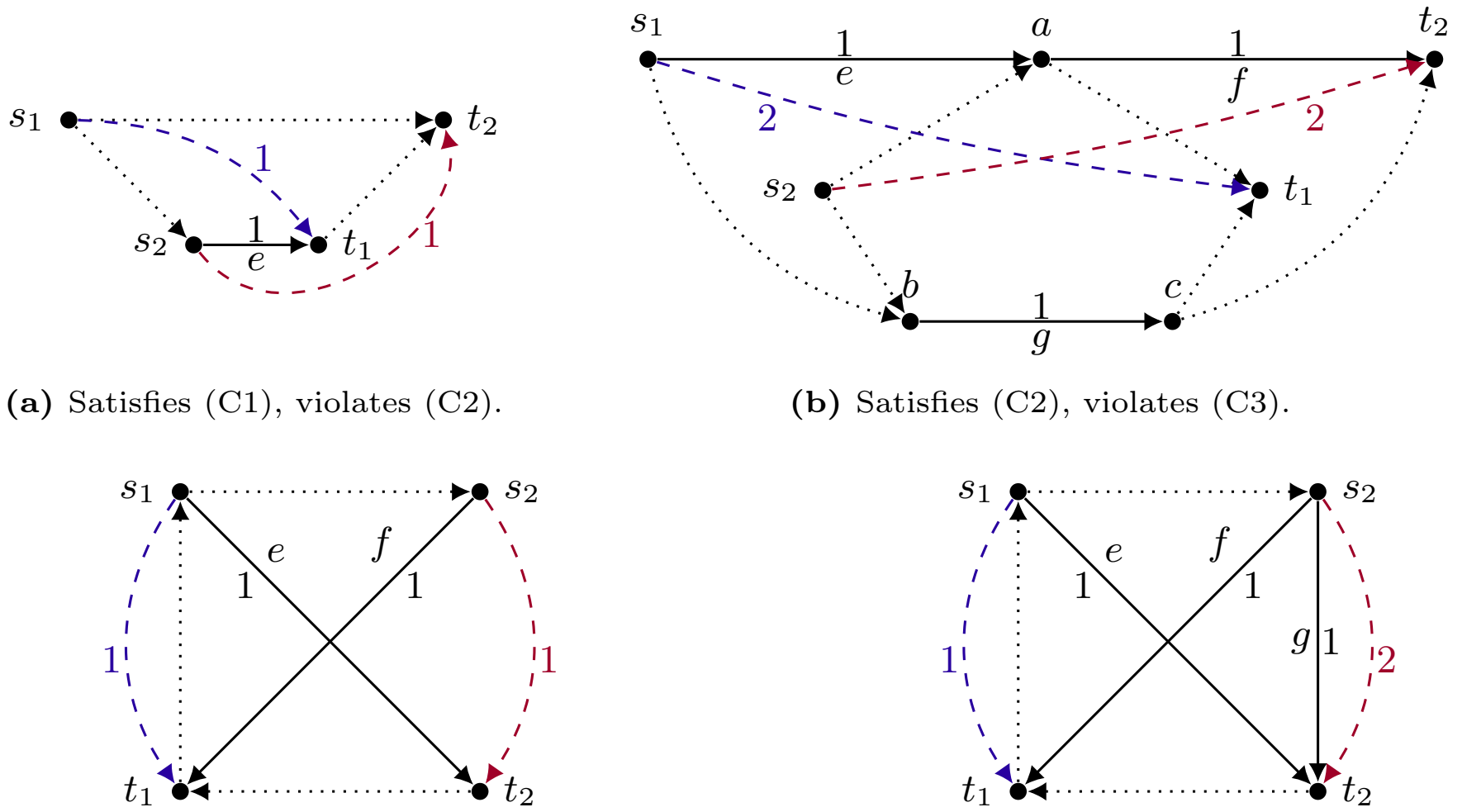
\begin{figure}[!t]
    \centering
    \begin{subfigure}[b]{0.35\textwidth}
        \centering
        \begin{tikzpicture}[scale=.8]
          \begin{scope}[scale=1.25]
            \node[solidnode, label=left:{$s_1$}] (s) at (0, 0) {};
            \node[solidnode, label=left:{$s_2$}] (v) at (1, -1) {};
            \node[solidnode, label=right:{$t_1$}] (w) at (2, -1) {};
            \node[solidnode, label=right:{$t_2$}] (t) at (3, 0) {};
          \end{scope}

            % \node[left=0ex of s, ccolor1] {\commodityinfo{s_1}{1}};
            % \node[below=0ex of v, ccolor2] {\commodityinfo{s_2}{1}};

            \draw
                (s) edge[dedge, dotted] (v)
                (v) edge[dedge] node[midway, above, inner sep=.1ex] {\footnotesize$1$} node[midway, below, inner sep=.4ex] {\footnotesize$e$} (w)
                (w) edge[dedge, dotted] (t)
                (s) edge[dedge, dotted] (t);

          \draw[demandedge,ccolor1, bend left=25] (s) edge node[pos=.75, above] {\footnotesize$1$} (w);
          \draw[demandedge,ccolor2] (v) .. controls ++(.9, -1.15) and ++(.25, -1.1) .. (t) node[pos=.65, right] {\footnotesize$1$};
        \end{tikzpicture}

        \caption{Satisfies~\eqref{C1}, violates~\eqref{C2}.}
        \label{hier:fig-c1}
    \end{subfigure}
    \hfill
    \begin{subfigure}[b]{0.6\textwidth}
        \centering
        \begin{tikzpicture}[scale=.70]
            \node[solidnode, label=above:{\footnotesize$s_1$}] (s1) at (0, 0) {};
            \node[solidnode, label=above:{\footnotesize$a$} ] (v) at (4.5, 0) {};
            \node[solidnode, label=above:{\footnotesize$t_2$}] (tt2) at (9, 0) {};
            \node[solidnode, label=left:{\footnotesize$s_2$}] (s2) at (2, -1.5) {};
            \node[solidnode, label=right:{\footnotesize$t_1$}] (tt1) at (7, -1.5) {};
            \node[solidnode,label=above:{\footnotesize$b$}  ] (v1) at (3, -3) {};
            \node[solidnode, label=above:{\footnotesize$c$} ] (w1) at (6, -3) {};

            % \node[above=0ex of s1, ccolor1] {\commodityinfo{s_1}{2}};
            % \node[above=-.25ex of s2, xshift=-2.2ex, ccolor2] {\commodityinfo{s_2}{2}};

            \draw
                (s1) edge[dedge] node[midway, above, inner sep=.1ex] {\footnotesize$1$} node[midway, below, inner sep=.4ex] {\footnotesize$e$} (v)
                (v) edge[dedge] node[midway, above, inner sep=.1ex] {\footnotesize$1$} node[midway, below, inner sep=.4ex] {\footnotesize$f$} (tt2)
                (s2) edge[dedge, dotted] (v)
                (v) edge[dedge, dotted] (tt1)
                (s2) edge[dedge, dotted] (v1)
                (w1) edge[dedge, dotted] (tt1)
                (s1) edge[dedge, dotted, bend right=30] (v1)
                (w1) edge[dedge, dotted, bend right=30] (tt2)
                (v1) edge[dedge] node[midway, above, inner sep=.1ex] {\footnotesize$1$} node[midway, below, inner sep=.4ex] {\footnotesize$g$} (w1);

            \draw[demandedge,ccolor1, bend right=5] (s1) edge node[pos=.2, below left, inner sep=.25ex] {\footnotesize$2$} (tt1);
            \draw[demandedge,ccolor2, bend right=5] (s2)  edge node[pos=.8, below right, inner sep=.25ex] {\footnotesize$2$} (tt2);
        \end{tikzpicture}
        \caption{Satisfies~\eqref{C2}, violates~\eqref{C3}.}
        \label{hier:fig-c2}
    \end{subfigure}
    \vspace{1em}
    
    \begin{subfigure}[b]{0.45\textwidth}
        \centering
        %%%%%%% BEGIN ALTERNATIVE FIGURE FOR C3-C4-Gadget %%%%%%%
        \begin{tikzpicture}[scale=.8]
            % \useasboundingbox (5,-1.5ex) rectangle ++(.1, -3);
            % \fill[green!10] (-.75, .75) rectangle (3.75, -3.75);
            \node[solidnode, label=left:{\footnotesize$s_1$}] (s1) at (0, 0) {};
            \node[solidnode, label=right:{\footnotesize$s_2$}] (s2) at (3, 0) {};
            \node[solidnode, label=right:{\footnotesize$t_2$}] (t2) at (3, -3) {};
            \node[solidnode, label=left:{\footnotesize$t_1$}] (t1) at (0, -3) {};

            \draw
                (s1) edge[dedge, dotted] (s2)
                (t1) edge[dedge, dotted] (s1)
                (t2) edge[dedge, dotted] (t1)
                (s1) edge[dedge] node[pos=.25, below left, inner sep=.1ex]{\footnotesize$1$} node[pos=.25, above right, inner sep=.4ex]{\footnotesize$e$} (t2)
                (s2) edge[dedge] node[pos=.25, below right, inner sep=.1ex]{\footnotesize$1$} node[pos=.25, above left, inner sep=.4ex]{\footnotesize$f$} (t1);

            \draw[demandedge,ccolor1, bend right=35] (s1) edge node[pos=.5, below left, inner sep=.25ex] {\footnotesize$1$} (t1);
            \draw[demandedge,ccolor2, bend left=35] (s2)  edge node[pos=.5, below right, inner sep=.25ex] {\footnotesize$1$} (t2);
        \end{tikzpicture}
        
        \caption{Satisfies~\eqref{C3}, violates~\eqref{C4}.}
        \label{hier:fig-c3}
    \end{subfigure}
    \hfill
    \begin{subfigure}[b]{0.45\textwidth}
        \centering
        
        \begin{tikzpicture}[scale=.8]
           
            \node[solidnode, label=left:{\footnotesize$s_1$}] (s1) at (0, 0) {};
            \node[solidnode, label=right:{\footnotesize$s_2$}] (s2) at (3, 0) {};
            \node[solidnode, label=left:{\footnotesize$t_1$}] (t1) at (0, -3) {};
            \node[solidnode, label=right:{\footnotesize$t_2$}] (t2) at (3, -3) {};

            \draw
                (s1) edge[dedge, dotted] (s2)
                (t1) edge[dedge, dotted] (s1)
                (t2) edge[dedge, dotted] (t1)
                (s1) edge[dedge] node[pos=.25, below left, inner sep=.1ex]{\footnotesize$1$} node[pos=.25, above right, inner sep=.4ex]{\footnotesize$e$} (t2)
                (s2) edge[dedge] node[pos=.25, below right, inner sep=.1ex]{\footnotesize$1$} node[pos=.25, above left, inner sep=.4ex]{\footnotesize$f$} (t1)
                (s2) edge[dedge] node[pos=.5, right, inner sep=.1ex]{\footnotesize$1$} node[pos=.5, left, inner sep=.3ex]{\footnotesize$g$} (t2);

            \draw[demandedge,ccolor1, bend right=35] (s1) edge node[pos=.5, below left, inner sep=.25ex] {\footnotesize$1$} (t1);
            \draw[demandedge,ccolor2, bend left=35] (s2)  edge node[pos=.5, below right, inner sep=.25ex] {\footnotesize$2$} (t2);
        \end{tikzpicture}
       
        \caption{Satisfies~\eqref{C4}, violates~\eqref{C5}.}
        \label{hier:fig-c4}
    \end{subfigure}

    \caption{Infeasible instances separating successive levels of the hierarchy.  Solid edges have capacity~$1$ and dotted edges capacity~$d_1+d_2$; dashed edges represent commodities with the demand given next to the respective edges.  Panels~(a) and (b) are adapted from~\cite{majthoub2026integer}.}
\label{hier:fig-examples}
\end{figure}

\begin{theorem}\label{hier:strict}
Even when restricted to two-commodity topologies,
\[
 \mathcal S_1\subsetneq\mathcal S_2\subsetneq\mathcal S_3
 \subsetneq\mathcal S_4\subsetneq\mathcal S_5.
\]
More precisely, for $r=1,2,3,4$, the~$r$th topology in~\cref{hier:fig-examples} belongs to~$\mathcal S_{r+1}\setminus\mathcal S_r$.
The displayed weights satisfy the~$r$th condition and violate the~$(r+1)$st.
\end{theorem}

The proof of \cref{hier:strict} can be found in \cref{hier:appendix-families} together
with sufficient families for the intermediate levels. In particular,
unique commodity paths ensure~\eqref{C2}-sufficiency, as the next lemma
shows; its proof is deferred to \cref{app:unique-paths}.

\begin{lemma}[Unique commodity paths]\label{hier:unique}
If every commodity has at most one simple directed source--sink path,
then $(G,H)\in\mathcal S_2$.
\end{lemma}

Checking~\eqref{C1} can be done in polynomial-time for a fixed number of
commodities~\cite{Poremba_Shepherd_2025}; checking~\eqref{C2}-\eqref{C4}
is strongly coNP-complete even for two commodities, while \eqref{C5}
is polynomial-time checkable (\cref{hier:easy-checks,hier:hard-checks}).
Recognizing the sufficient topology classes, however, is a separate problem. Poremba and Shepherd \cite[Thm.~8]{Poremba_Shepherd_2025} proved that
recognizing~$\mathcal S_1$ is coNP-hard already for roundtrip demands.
We show hardness at every level (see~\cref{hier:appendix-recognition} for proof).

\begin{theorem}[Recognition]\label{hier:recognition}
For each $r\in\{1,2,3,4\}$, deciding whether a topology belongs
to~$\mathcal S_r$ is in coNP, for any number of commodities. It is
coNP-hard already for two-commodity topologies with two-matching demands;
more strongly, it is coNP-hard to distinguish topologies
in~$\mathcal S_1$ from topologies outside~$\mathcal S_4$.
\end{theorem}

\section{From length functions to forbidden relevant minors}
\label{char:section}

We now characterize the two-commodity topologies in~$\mathcal{S}_1$, completing the work of
Poremba and Shepherd~\cite{Poremba_Shepherd_2025}.

\begin{theorem}\label{char:main}
A directed two-commodity topology is cut-sufficient if and only if it
contains neither the bad triangle nor the bad dual triangles as a
relevant minor (\cref{char:relevant}).
\end{theorem}

\subsection{Relevant minors and the known two-path case}
\label{char:sec-minors}

A \emph{minor} of~$(G,H)$ is obtained by deleting supply and demand
edges and contracting a set~$F$ of supply edges, identifying their
endpoints in both~$G$ and~$H$. We treat~$F$ and the deleted edges as
part of the minor, since the following extension refers to them.

\begin{definition}[Extension of weights]\label{char:extension}
Let~$(G',H')$ be a minor of~$(G,H)$ and let~$(u',d')$ be weights
on~$(G',H')$.  Let $D:=\sum_{i}d'_i$.  The \emph{extension}~$(u^{\mathrm{ext}},d^{\mathrm{ext}})$ 
of~$(u',d')$ to~$(G,H)$ is
\[
 u^{\mathrm{ext}}_e:=
 \begin{cases}
  0,& e\ \text{deleted},\\
  D,& e\in F,\\
  u'_e,& \text{otherwise},
 \end{cases}
 \qquad
 d^{\mathrm{ext}}_i:=
 \begin{cases}
  0,& i\ \text{deleted},\\
  d'_i,& \text{otherwise}.
 \end{cases}
\]
\end{definition}

\begin{definition}[Relevant minor,
 {\cite[Def.~2]{Poremba_Shepherd_2025}}]\label{char:relevant}
A minor is \emph{relevant} if the extension of every weighting of the
minor that satisfies~\eqref{C1} also satisfies~\eqref{C1} on~$(G,H)$.
\end{definition}
Poremba and Shepherd use capacity~$+\infty$ on contracted edges.
Capacity~$D$ gives the same relevant minors: any cut containing such an
edge has capacity at least the total demand. Our convention keeps the
weights finite and integral. As in~\cite{Poremba_Shepherd_2025},
relevance and cut-sufficiency are unchanged if weights are restricted
to integers.\footnote{A failure of relevance is described by homogeneous rational
cut inequalities and one strict reverse cut inequality, and a failure of
cut-sufficiency analogously by one violated integral length inequality
from~\eqref{C5}.  A nonempty system of either kind has a rational
solution, which can be scaled to an integral one.}
Deletions always give relevant minors, contractions need not~\cite[Fig.~4]{Poremba_Shepherd_2025}, and relevance is transitive because extensions compose. A routing on the extended
weights in the original topology projects to a routing on the minor, so a
relevant minor outside~$\mathcal S_1$ certifies that its parent is
outside~$\mathcal S_1$.

Demand edges~$(a,b)$, $(b,c)$
are called \emph{two-path demands} when~$a,b,c$ are distinct, and
\emph{roundtrip demands} when $c=a$.
We use the following results of Poremba and
Shepherd~\cite{Poremba_Shepherd_2025}.

\begin{theorem}[{\cite[Thms.~2, 3, and~5]{Poremba_Shepherd_2025}}]\label{char:ps}
For two-path or roundtrip demands, the following are equivalent:
\begin{enumerate}
\renewcommand{\labelenumi}{(\roman{enumi})}
\item the topology is cut-sufficient;
\item every simple path for one commodity is edge-disjoint from every
      simple path for the other;
\item neither the bad triangle nor the bad dual triangles is a
      relevant minor.
\end{enumerate}
\end{theorem}

Thus, it only remains to consider \emph{two-matching demands} with four
distinct terminals. Every other two-commodity demand graph has two-path or
roundtrip demands, or a common source or sink.  In the last case, the topology is cut-sufficient by the single-commodity max-flow/min-cut theorem applied after adding a super-sink or a super-source. For any deceptive weights on a two-matching demands topology, both demands are positive, and each~$t_i$ is reachable from~$s_i$ by~\eqref{C1}.

\begin{lemma}[Cross-commodity contraction]\label{char:safe}
Let~$(G,H)$ have two-matching demands, and let~$P$ be a directed~$s_i$--$t_j$~path, $i\ne j$. Contracting~$P$ produces a relevant minor.
If~$P$ avoids~$t_i$ and~$s_j$, it produces two-path demands~$s_j\to z\to t_i$, where~$z$ is the contracted vertex.
\end{lemma}
\begin{proof}
Let weights on the contracted topology satisfy~\eqref{C1}, and let
$Y\subseteq V$. If~$Y$ contains all or none of~$V(P)$, its cut inequality
for the extension is that of the corresponding set in the minor. If an
edge of~$P$ leaves~$Y$, its capacity~$D$ covers all demand. Otherwise~$P$
only enters~$Y$ (and does not leave it), so $s_i\notin Y$ and $t_j\in Y$, thus no demand
leaves~$Y$.
\qed
\end{proof}

\subsection{A short cross-commodity path} \label{char:sec-cross}

Fix a two-matching instance~$(G,H,u,d)$ with deceptive weights. Since it
is infeasible, \eqref{C5} fails: there is a \emph{witness}, an integral
length vector~$\ell\ge0$ with

\begin{equation}\label{char:witness}
 \sum_e u_e\ell_e<d_1m_1+d_2m_2,
 \qquad m_i:=\dist_\ell(s_i,t_i),
\end{equation}
where each~$m_i$ is finite because~$t_i$ is reachable from~$s_i$. In
\cref{char:sec-layer} we choose a witness of minimum total length; the
next lemma holds for every witness.

\begin{lemma}\label{char:short-cross}
Every witness satisfies
$
 \dist_\ell(s_1,t_2)<m_1
 \quad\text{or}\quad
 \dist_\ell(s_2,t_1)<m_2.
$
\end{lemma}

\begin{proof}
Suppose both inequalities fail.  Set~$\bar m:=\max\{m_1,m_2\}$ and define 
\[
 \sigma(v):=\min\{\dist_\ell(v,t_1),\dist_\ell(v,t_2),\bar m\}.
\]
Then $\sigma(t_1)=\sigma(t_2)=0$ and $\sigma(s_i)=m_i$ for $i=1,2$.  For
each integer~$0\le\alpha<\bar m$, the set~$Y_\alpha:=\{v:\sigma(v)>\alpha\}$
contains no sink, and it contains~$s_i$ exactly when $\alpha<m_i$.
Summing the \eqref{C1} inequalities of $\delta_G^+(Y_\alpha)$ gives
\begin{equation} \label{eq:sum Y_alpha}
 \sum_{\alpha=0}^{\bar m-1}u(\delta_G^+(Y_\alpha))\ge d_1m_1+d_2m_2.
\end{equation}
Every edge~$(v,w)$ satisfies $\sigma(v)\le\ell_{vw}+\sigma(w)$, since
each term of the minimum does. Hence~$(v,w)$ leaves at most~$\ell_{vw}$
of the sets~$Y_\alpha$, and the left-hand side of~\eqref{eq:sum Y_alpha}
is at most $\sum_eu_e\ell_e$, a contradiction.
\qed
\end{proof}

\subsection{Every layer interacts with both commodities}
\label{char:sec-layer}

Choose a witness~$\ell$ of minimum total length~$\|\ell\|_1$. For every
commodity~$i$, let $\rho_i(v):=\min\{\dist_\ell(s_i,v),m_i\}$; it is
finite everywhere and $\rho_i(w)-\rho_i(v)\le\ell_{vw}$ on every
edge. For an integer $0\le\alpha<m_i$, define
$
 X_i^\alpha :=\{v:\rho_i(v)\le\alpha\}$ and $
 L_i^\alpha :=\{(x,y)\in\delta_G^+(X_i^\alpha):
        \text{$y$ reaches~$t_i$ in~$G[V\setminus X_i^\alpha]$}\}.
$
We call~$L_i^\alpha$ a \emph{layer} of commodity~$i$: it consists of the
edges through which an $s_i$--$t_i$ path can leave~$X_i^\alpha$ for the
last time. Every $s_i$--$t_i$ path has a last exit, so~$L_i^\alpha$
blocks commodity~$i$ and contains a cut separating it
(\cref{sec:preliminaries}). Hence~\eqref{C1} gives
\begin{equation}\label{char:layer-capacity}
 u(L_i^\alpha)\ge d_i.
\end{equation}
Each edge $(x,y)\in L_i^\alpha$ has $\rho_i(y)\ge\alpha+1>\rho_i(x)$,
hence $\ell_{xy}\ge1$. It also lies on a simple commodity-$i$ path: a
shortest $s_i$--$x$ path inside~$X_i^\alpha$, then~$(x,y)$, then a simple
$y$--$t_i$ path outside~$X_i^\alpha$.

\begin{lemma}[No private layer]\label{char:layer}
Let~$\ell$ be a witness of minimum total length, let $i\ne j$, and let
$0\le\alpha<m_i$ be an integer.  Put $\ell':=\ell-\ind{L_i^\alpha}$ and
$m'_h:=\dist_{\ell'}(s_h,t_h)$ for $h=1,2$.  Then,
$m'_i=m_i-1$ and $m'_j<m_j$.
\end{lemma}
\begin{proof}
The vector~$\ell'$ is nonnegative and integral.  An $\ell$-shortest
commodity-$i$ path uses an edge of~$L_i^\alpha$, so $m'_i\le m_i-1$.

For the converse, let $\phi(v):=\rho_i(v)$ if $v\in X_i^\alpha$ and
$\phi(v):=\rho_i(v)-1$ otherwise. Then $\phi(w)-\phi(v)\le\ell'_{vw}$ for
every edge~$(v,w)$. If both ends lie on the same side of~$X_i^\alpha$,
this is the inequality for~$\rho_i$, since $\ell'=\ell$ there. If $(v,w)$
leaves~$X_i^\alpha$, then
$\phi(w)-\phi(v)=\rho_i(w)-\rho_i(v)-1\le\ell_{vw}-1\le\ell'_{vw}$. If it
enters~$X_i^\alpha$, then $\phi(w)\le\alpha<\phi(v)+1$, so
$\phi(w)-\phi(v)\le0$ as $\phi$ is integral. Since $\phi(s_i)=0$ and $\phi(t_i)=m_i-1$, summing
along any $s_i$--$t_i$ path gives $m'_i\ge m_i-1$. 

The layer~$L_i^\alpha$ is nonempty, so $\|\ell'\|_1<\|\ell\|_1$.  Minimality of~$\ell$ implies that~$\ell'$ is not a witness.  Using~\eqref{char:layer-capacity}, we obtain
\[
\begin{aligned}
 \sum_e u_e\ell_e
 &=\sum_e u_e\ell'_e+u(L_i^\alpha)
 \ge d_i(m_i-1)+d_jm'_j+d_i
 =d_im_i+d_jm'_j.
\end{aligned}
\]
Comparison with~\eqref{char:witness}, and $d_j>0$, yields $m'_j<m_j$.
\qed
\end{proof}

Thus a shortest commodity-$j$ path for~$\ell'$ uses the layer.
We next show that a suitable layer supplies a common edge that survives
the required contraction.

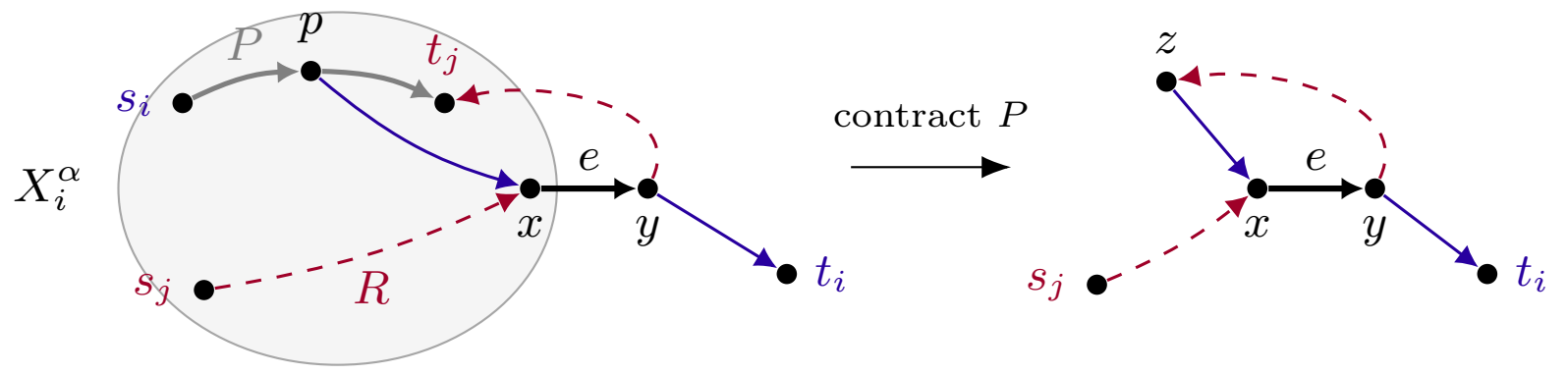
\begin{figure}[tb]
\centering
\begin{tikzpicture}[x=.72cm,y=.72cm,
  pathp/.style={dedge,black!50,line width=.25ex},
  pathi/.style={dedge,ccolor1},
  pathr/.style={dedge,ccolor2,dashed}]
  \useasboundingbox (-4,1.5) rectangle (10,-1.5);
  % \fill[red!10] (-4,1.5) rectangle (10,-1.5);
  \draw[fill=black!4,draw=black!35] (-1.55,0) ellipse (2.05 and 1.65);
  \node[font=\small] at (-4.25,0) {$X_i^\alpha$};
  \node[solidnode,label=left:{\color{ccolor1}$s_i$}] (si) at (-3,.8) {};
  \node[solidnode,label=above:$p$] (p) at (-1.8,1.1) {};
  \node[solidnode,label=above:{\color{ccolor2}$t_j$}] (tj) at (-.55,.8) {};
  \node[solidnode,label=left:{\color{ccolor2}$s_j$}] (sj) at (-2.8,-.95) {};
  \node[solidnode,label=below:$x$] (x) at (.25,0) {};
  \node[solidnode,label=below:$y$] (y) at (1.35,0) {};
  \node[solidnode,label=right:{\color{ccolor1}$t_i$}] (ti) at (2.65,-.8) {};
  \draw[pathp] (si) to[bend left=12] node[above,font=\small] {$P$} (p);
  \draw[pathp] (p) to[bend left=12] (tj);
  \draw[pathi] (p) to[bend right=12] (x);
  \draw[pathr] (sj) to[bend right=8] node[below,font=\small] {$R$} (x);
  \draw[dedge,line width=.3ex] (x) -- node[above] {$e$} (y);
  \draw[pathi] (y) -- (ti);
  \draw[pathr] (y) to[out=65,in=15] (tj);
  \draw[-{Latex[length=2mm]}] (3.25,.2) -- (4.75,.2);
  \node[font=\scriptsize,align=center] at (4,.7) {contract~$P$};
  \begin{scope}[shift={(6.2,0)}]
  \node[solidnode,label=above:$z$] (z2) at (0,1) {};
  \node[solidnode,label=left:{\color{ccolor2}$s_j$}] (sj2) at (-.65,-.9) {};
  \node[solidnode,label=below:$x$] (x2) at (.85,0) {};
  \node[solidnode,label=below:$y$] (y2) at (1.95,0) {};
  \node[solidnode,label=right:{\color{ccolor1}$t_i$}] (ti2) at (3,-.8) {};
  \draw[pathi] (z2) -- (x2);
  \draw[pathr] (sj2) to[bend right=8] (x2);
  \draw[dedge,line width=.3ex] (x2) -- node[above] {$e$} (y2);
  \draw[pathi] (y2) -- (ti2);
  \draw[pathr] (y2) to[out=65,in=10] (z2);
  \end{scope}
\end{tikzpicture}
\caption{The contraction in~\cref{char:preserve}, shown schematically.
The cross-commodity path~$P$ (gray) lies inside~$X:=X_i^\alpha$; the first edge~$e$ of the commodity-$j$ path~$R$ (dashed) in the layer $L:=L_i^\alpha$ leaves~$X$.  The
retained commodity-$i$ path starts at its last contact~$p$ with~$P$ and
stays outside~$P$ thereafter.  The retained prefix of~$R$ ends at its
first contact with~$P$.  Both become simple paths through~$e$ after
contraction.}
\label{char:layer-figure}
\end{figure}

\begin{lemma}[Preserving a common edge on simple paths]\label{char:preserve}
Let~$\ell$ be a witness of minimum total length and suppose that, for
$i\ne j$, some~$\ell$-shortest $s_i$--$t_j$~path~$P$ has length~$\alpha<m_i$.  Then contracting~$P$ produces a relevant minor with
two-path demands and a common edge on simple paths for its two
commodities.
\end{lemma}
\begin{proof}
Define $X:=X_i^\alpha$, $L:=L_i^\alpha$, and $\ell':=\ell-\ind{L}$;
see~\cref{char:layer-figure}.
Since~$P$ is $\ell$-shortest of length~$\alpha$, every vertex of~$P$ has
$\rho_i\le\alpha$ and so belongs to~$X$; as $L\subseteq\delta^+_G(X)$,
the vectors $\ell$ and~$\ell'$ agree on~$P$.
Let~$R$ be a simple~$\ell'$-shortest~$s_j$--$t_j$~path.  By~\cref{char:layer}, $\ell'(R)<m_j$, so~$R$ uses an edge of~$L$; let
$e=(x,y)$ be its first.  In particular, $\ell$ and~$\ell'$ agree on~$R[s_j,x]$.

\emph{Claim: $R[s_j,x]$ avoids~$V(P)$.}  Otherwise, choose a vertex~$q$
in their intersection and write $\gamma:=\ell(P[s_i,q])=\rho_i(q)$,
where equality holds because~$P$ is shortest and prefixes of $\ell$-shortest paths are $\ell$-shortest.  Since
$\rho_i(y)\ge\alpha+1$,
\begin{align*}
 \ell'(R[q,y])
 &=\ell(R[q,x])+\ell_e-1\\
 &\ge (\rho_i(x)-\gamma)+ (\rho_i(y)-\rho_i(x))-1
 \ge \alpha-\gamma
 =\ell(P[q,t_j]).
\end{align*}
Replacing the suffix of~$R$ after~$q$ by~$P[q,t_j]$, and using that the
remaining~$\ell'$-length of~$R[y,t_j]$ is nonnegative, therefore gives an~$s_j$--$t_j$~walk of~$\ell$-length at most $\ell'(R)<m_j$.  Removing
cycles contradicts the definition of~$m_j$.  This proves the claim, and
in particular $s_j\notin V(P)$.  Also $t_i\notin V(P)$, since
$\ell(P)<m_i$.

We may therefore contract~$P$ to produce a relevant minor by \cref{char:safe}.
The edge~$e$ survives: its tail is outside~$V(P)$ by the claim, and
its head lies outside~$X\supseteq V(P)$.  To obtain a simple
commodity-$j$ path in the minor, notice that since $t_j\in V(P)$, the path~$R$ meets~$V(P)$; retain~$R$ only up to its
first such vertex.  This prefix contains~$e$ and becomes a simple~$s_j$--$z$~path after contraction.

For commodity~$i$, take a simple shortest~$s_i$--$x$~path~$Q$ within~$X$ and a simple $y$--$t_i$~path~$T$ outside~$X$.
Let~$p$ be the last vertex of~$Q$ in~$V(P)$.  The path~$Q[p,x]\cdot e\cdot T$ meets~$V(P)$ only at~$p$ and is simple, because~$Q$ and~$T$ lie on opposite sides of~$X$.  After contraction it is a
simple~$z$--$t_i$~path containing~$e$, as required.
\qed
\end{proof}

\begin{proof}[of~\cref{char:main}]
Both forbidden topologies admit deceptive weights, so
necessity follows from the closure of cut-sufficiency under relevant
minors.  For sufficiency, suppose that~$(G,H)$ is not cut-sufficient and
choose deceptive weights.  \Cref{char:ps} handles
two-path and roundtrip demands, and topologies with a common source or sink are cut-sufficient and hence excluded, so the demands are two-matching and positive.  Choose a witness~$\ell$ of minimum total length.  By~\cref{char:short-cross}, there are
$i\ne j$ with $\alpha:=\dist_\ell(s_i,t_j)<m_i$; let~$P$ be an~$\ell$-shortest~$s_i$--$t_j$~path.  \Cref{char:preserve} implies that contracting~$P$ yields a relevant minor with two-path demands and a common edge on
simple paths for both commodities.  By~\cref{char:ps}, this relevant minor
contains a forbidden relevant minor, and by transitivity so does~$(G,H)$.
\qed
\end{proof}

\section{Integrality at the first four levels}
 \label{int:section}

Nagamochi and Ibaraki~\cite{Nagamochi_Ibaraki_89} proved that every
cut-sufficient topology is integrally cut-sufficient.  We extend this to
the next three levels of the hierarchy and make it constructive.
For~\eqref{C2} and~\eqref{C3}, the existence statement was already
observed in~\cite{majthoub2026integer} without proof. Since
$\mathcal S_1\subseteq\dots\subseteq\mathcal S_4$, it suffices to
treat~$\mathcal S_4$. 

\begin{theorem}\label{int:main}
%Let $r\in\{1,2,3,4\}$ and
Let $(G,H)\in\mathcal S_4$. 

Every feasible integral instance
on~$(G, H)$ admits an integral multiflow, and one can be computed in
polynomial time.
\end{theorem}

On~$\mathcal S_r$ with $r\le4$, feasibility is equivalent
to~\textup{(C$r$)}, so the theorem states that sufficiency
of~\textup{(C$r$)} implies integral sufficiency. The algorithm assumes the
topology belongs to~$\mathcal S_r$; it does not recognize that class. The algorithm only tests fractional feasibility and
solves multiflow linear programs; it never checks~\textup{(C$r$)},
which is coNP-hard for $r\ge2$.

For a multiflow~$x$, let $x_e:=\sum_i x^i_e$ be the total flow value on edge~$e$, and for $F\subseteq E$ let $x(F):=\sum_{e\in F}x_e$. 
%The property of~\eqref{C4} that drives this section is that all its capacity coefficients are~$0$
% or~$1$. 
The first four cut conditions share the following common feature that drives this section: 

\begin{lemma}[Unit reduction]\label{int:unit}
Let $r\in\{1,2,3,4\}$, let~$(u,d)$ be integral, let~$x$ be a feasible multiflow
for~$(G,H,u,d)$, and let~$e$ be an edge with $x_e<u_e$.  Then $(u-\ind{\{e\}},d)$ satisfies~\textup{(C$r$)}.
\end{lemma}

The proof of \cref{int:unit} is deferred to \cref{apx:missing-proofs}. It uses that every capacity coefficient is $0$ or~$1$, which fails
for~\eqref{C5}.  So does the conclusion of~\cref{int:main}:
$\mathcal S_5$ contains every topology, whereas fractional feasibility
does not imply integral feasibility, even with two commodities.

The next lemma makes the total flow on every edge integral, although
individual commodity flows may remain fractional. For fixed demands, it
extends \cite[Lem.~3.3]{Nagamochi_Ibaraki_89}.
Its proof is also deferred to \cref{apx:missing-proofs}.

\begin{lemma}[Integral total flow]\label{int:total}
%Let $r\in\{1,2,3,4\}$ and 
Fix an integral demand vector~$d$, and suppose that $(G,H,u',d)$ is feasible for every integral capacity
vector~$u'$ such that $(u',d)$ satisfies~\eqref{C4}.
Then, for every integral~$u$ such that $(u,d)$ is (fractionally) feasible, one can compute, using $|E|$ linear programs,
an integral capacity vector~$\tilde u\le u$ such that
$(G,H,\tilde u,d)$ is feasible and every feasible multiflow
under~$\tilde u$ saturates every edge.
\end{lemma}

\begin{proof}[of~\cref{int:main}]
\emph{Correctness.} We construct an integral multiflow by induction on
the total demand; the case $d=0$ is trivial. As
$(G,H)\in\mathcal S_4$, the hypothesis of~\cref{int:total} holds for
every integral~$d$, so the lemma yields integral $\tilde u\le u$ and a
feasible multiflow~$x$ with $x_e=\tilde u_e$ for every edge~$e$. Flow
conservation of~$x$ gives, for every vertex~$v$,
\begin{equation}\label{int:euler}
 \tilde u(\delta_G^+(v))-\tilde u(\delta_G^-(v))
 =\sum_{i:\,s_i=v}d_i-\sum_{i:\,t_i=v}d_i.
\end{equation}
Choose~$i$ with $d_i>0$. Removing one unit of commodity~$i$ from~$x$
shows that $(\tilde u,d-\ind{\{i\}})$ is feasible, so by induction it
admits an integral multiflow~$z$. The unused capacities
$\bar u:=\tilde u-\sum_{i' \in [k]}z^{i'}$ are nonnegative integers, and
subtracting the conservation equations of~$z$ from~\eqref{int:euler}
shows that~$\bar u$ has net outflow~$1$ at~$s_i$, $-1$ at~$t_i$,
and~$0$ elsewhere. Hence~$t_i$ is reachable from~$s_i$ along edges satisfying
$\bar u_e\ge1$: otherwise the set of vertices so reachable would have
net outflow~$1$ but no outgoing unused capacity. Adding one unit of
commodity~$i$ along such a path to~$z$ gives an integral multiflow
for~$(\tilde u,d)$, and hence for~$(u,d)$.

\emph{Running time.} The recursion makes one call per unit of demand,
so we first reduce the demands. Compute a feasible multiflow~$x$,
cancel its cycles, and decompose each commodity into at most $|E|$
paths. Rounding these path flows down gives an integral multiflow~$y$.
The residual instance with $d^{\mathrm{res}}_i:=d_i-\sum_Py^i_P$ and
$u^{\mathrm{res}}_e:=u_e-\sum_iy^i_e$ is integral and feasible,
witnessed by $x-y$, and $d^{\mathrm{res}}_i< |E|$ because it is a sum of
at most~$|E|$ fractional parts. Apply the recursion to it and add~$y$.
Fewer than~$k|E|$ recursive steps remain, each solving~$|E|$ linear
programs and finding one path, so at most $1+k|E|^2$ linear programs are
solved in total. Their encoding sizes stay polynomial, since the
topology is unchanged and capacities and demands only decrease.
\qed
\end{proof}

In particular, \cref{int:main} implies that integral multiflows can be computed in polynomial time,
whenever they exist, on parallel-channel topologies, cactus and
block-ESP digraphs, and topologies in which all but one commodity have
a unique path (\cref{hier:parallel,hier:cactus,hier:all-but-one}).

\section*{Acknowledgements: AI Disclosure}

The main result of this paper, the proof of the Poremba--Shepherd
conjecture, was obtained without any use of AI, as was the generalization
of the Nagamochi--Ibaraki theorem to conditions~\eqref{C2}--\eqref{C4}.
We subsequently used AI tools to help simplify arguments or make them
more precise and, in particular, to improve the language and presentation.
Through iterative exchanges combining the authors' mathematical
expertise with AI-generated suggestions, AI also contributed to the
hardness results for recognizing sufficient topologies, extensions and
clearer formulations of sufficiency arguments, constructions for gaps at
the higher levels of the hierarchy, and results on length tests with
bounded coefficients. The authors take full responsibility for all
results, proofs, and other content of the paper.

\bibliographystyle{splncs04}
\bibliography{bibliography}

\appendix
\crefalias{section}{appendix}

\renewcommand{\theHsection}{appendix.\Alph{section}}
\section{Sufficient families and strictness of the hierarchy}
\label{hier:appendix-families}

This appendix proves the membership results stated in~\cref{hier:sec}
and, alongside them, verifies the four examples
of~\cref{hier:fig-examples}, which proves~\cref{hier:strict}. In these
examples every dotted edge has capacity $D:=d_1+d_2$. Hence an edge set
containing a dotted edge satisfies every inequality
of~\eqref{C1}--\eqref{C3}. Lowering the length of a dotted edge by one
lowers the left-hand side of~\eqref{C4} or~\eqref{C5} by~$D$ and each
commodity distance by at most one, since a shortest path uses the edge at most
once; it preserves any violation. Below we therefore check only edge
sets and length vectors that avoid the dotted edges.

\subsection{All or all but one commodity has a unique path}
\label{app:unique-paths}

\begin{proof}[of \cref{hier:unique}]
Let~$(u,d)$ satisfy~\eqref{C2}, and omit zero-demand commodities.
Every remaining commodity~$i$ has a path, since otherwise the
cut~$\delta^+_G(X_i(\emptyset))=\emptyset$ separates it and
violates~\eqref{C1}; let~$P_i$ be this path. Feasibility is
equivalent to
\begin{equation}\label{eq:capacity bound}
 u_e\ge\sum_{i:\,e\in P_i}d_i\qquad\text{for }e\in E.
\end{equation}
Fix an edge~$e=(v,w)$ on some~$P_i$. The graph~$G-e$ has no
$v$--$w$ path, since splicing one into~$P_i$ and shortcutting cycles
would give a second simple $s_i$--$t_i$ path. Hence the set~$X$ reachable
from~$v$ in~$G-e$ satisfies $\delta^+_G(X)=\{e\}$. This cut blocks
exactly the commodities whose path uses~$e$, so~\eqref{C2} for~$X$
is~\eqref{eq:capacity bound} for~$e$.
\qed
\end{proof}
\begin{corollary}\label{hier:ex1}
The first topology of~\cref{hier:fig-examples} lies in~$\mathcal S_2\setminus\mathcal S_1$.
\end{corollary}

\begin{proof}
By \cref{hier:unique}, since both commodities have unique paths, the topology lies
in~$\mathcal S_2$. For the given weights, a cut
that separates one commodity contains an edge of its path, and a cut
that separates both contains the dotted edge $(s_1,t_2)$; so~\eqref{C1}
holds. The cut $\delta^+_G(\{s_2\})=\{(s_2,t_1)\}$ blocks both
commodities and violates~\eqref{C2}.
\qed
\end{proof}

\begin{lemma}\label{hier:all-but-one}
If all but at most one commodity have at most one simple directed
source--sink path, then $(G,H)\in\mathcal S_4$.
\end{lemma}
\begin{proof}
Let~$(u,d)$ satisfy~\eqref{C4}, and omit zero-demand commodities; by~\eqref{C1} every 
positive-demand commodity has a path.
If all remaining commodities have a unique path, then \cref{hier:unique} applies.
 Otherwise, let~$j$ be the remaining commodity with
several paths, let~$I$ be the others, with unique paths~$P_i$, and let
$b_e:=\sum_{i\in I:\,e\in P_i}d_i$.
Applying~\eqref{C4} to the length vector~$\ind{\{e\}}$ gives
$u_e\ge b_e$, so the commodities in~$I$ can be routed on their paths.
For any $X\subseteq V$ with $s_j\in X$ and $t_j\notin X$, apply~\eqref{C4} to $\ell=\ind{\delta_G^+(X)}$ to obtain
\[
 u(\delta_G^+(X))
 \ge\sum_{i\in I}d_i|P_i\cap\delta_G^+(X)|+d_j
 =b(\delta_G^+(X))+d_j.
\]
Thus the residual capacities~$u-b$ satisfy every~$s_j$--$t_j$~cut
inequality.  The max-flow/min-cut theorem routes commodity~$j$ in
these residual capacities.
\qed
\end{proof}

\begin{corollary}\label{hier:ex3}
The third topology of~\cref{hier:fig-examples} lies
in~$\mathcal S_4\setminus\mathcal S_3$.
\end{corollary}
\begin{proof}
Commodity~$2$ has the unique path $s_2\to t_1\to s_1\to t_2$, so the
topology lies in~$\mathcal S_4$ by~\cref{hier:all-but-one}. Let
$e=(s_1,t_2)$ and $f=(s_2,t_1)$. Deleting~$e$ or~$f$ alone blocks only
commodity~$2$, and deleting both blocks both; so~\eqref{C3} holds for the given weights. The
lengths $\ell_e=\ell_f=1$ give distances~$1$ and~$2$ for commodities~$1$~and~$2$, respectively; therefore violating~\eqref{C4}.
\qed
\end{proof}

\begin{corollary}\label{hier:ex4}
The fourth topology of~\cref{hier:fig-examples} lies
in~$\mathcal S_5\setminus\mathcal S_4$.
\end{corollary}
\begin{proof}
Every topology lies in~$\mathcal S_5$. The fourth topology is the third topology with the added edge $g=(s_2,t_2)$. We first see that for the given weights~\eqref{C4} holds. Assuming $0$ lengths
on the dotted edges,
\[
 \dist_\ell(s_1,t_1)=\min\{\ell_e,\ell_f,\ell_g\},\qquad
 \dist_\ell(s_2,t_2)=\min\{\ell_g,\ell_e+\ell_f\}.
\]
For $\ell_e,\ell_f,\ell_g\in\{0,1\}$,
\begin{equation}\label{hier:eq-gadget}
 \ell_e+\ell_f+\ell_g\ge\min\{\ell_e,\ell_f,\ell_g\}
 +2\min\{\ell_g,\ell_e+\ell_f\}.
\end{equation}
\Cref{hier:eq-gadget} is immediate if $\ell_g=0$. If $\ell_g=1$, the cases
$\ell_e+\ell_f=0,1,2$ give right-hand sides $0,2,3$ (and left-hand sides $1,2,3$), respectively. Hence~\eqref{C4}
holds. The lengths $(\ell_e,\ell_f,\ell_g)=(1,1,2)$ give capacity
cost~$4$ and weighted distance $1+2\cdot2=5$, violating~\eqref{C5}.
\qed
\end{proof}

\subsection{Parallel channels} \label{families appendix: parallel channels}
\begin{figure}[!b]
\centering
\begin{tikzpicture}[x=1cm,y=.7cm,
  connector/.style={dedge,line width=.11ex},
  channel/.style={dedge,line width=.24ex}]
  \draw[rounded corners=5pt,fill=black!4,draw=black!20]
    (3.2,1.25) rectangle (6.8,1.85);
  \draw[rounded corners=5pt,fill=black!4,draw=black!20]
    (3.2,-.3) rectangle (6.8,.3);
  \draw[rounded corners=5pt,fill=black!4,draw=black!20]
    (4.7,-1.85) rectangle (5.3,-1.25);

  \node[solidnode,label=left:{\color{ccolor1}$s_1$}] (s1) at (0,1) {};
  \node[solidnode,label=left:{\color{ccolor2}$s_2$}] (s2) at (0,-1) {};
  \node[solidnode,label=right:{\color{ccolor1}$t_1$}] (t1) at (10,1) {};
  \node[solidnode,label=right:{\color{ccolor2}$t_2$}] (t2) at (10,-1) {};
  \node[solidnode] (a1) at (3.4,1.55) {};
  \node[solidnode] (b1) at (5,1.55) {};
  \node[solidnode] (c1) at (6.6,1.55) {};
  \node[solidnode] (a2) at (3.4,0) {};
  \node[solidnode] (c2) at (6.6,0) {};
  \node[solidnode] (a3) at (5,-1.55) {};

  \foreach \first in {a1,a2,a3} {
    \draw[connector,ccolor1] (s1) -- (\first);
    \draw[connector,ccolor2] (s2) -- (\first);
  }
  \foreach \last in {c1,c2,a3} {
    \draw[connector,ccolor1] (\last) -- (t1);
    \draw[connector,ccolor2] (\last) -- (t2);
  }
  \draw[channel] (a1) -- (b1);
  \draw[channel] (b1) -- (c1);
  \draw[channel] (a2) -- (c2);
  \node[above,font=\small] at (5,1.85) {$W_1$};
  \node[above,font=\small] at (5,.3) {$W_2$};
  \node[above,font=\small] at (5,-1.25) {$W_3$};
\end{tikzpicture}
\caption{A parallel-channel topology with two commodities and three
channels (shaded), including the single-vertex channel~$W_3$.
Colors distinguish the terminal pairs and their connectors.}
\label{hier:fig-parallel-channels}
\end{figure}
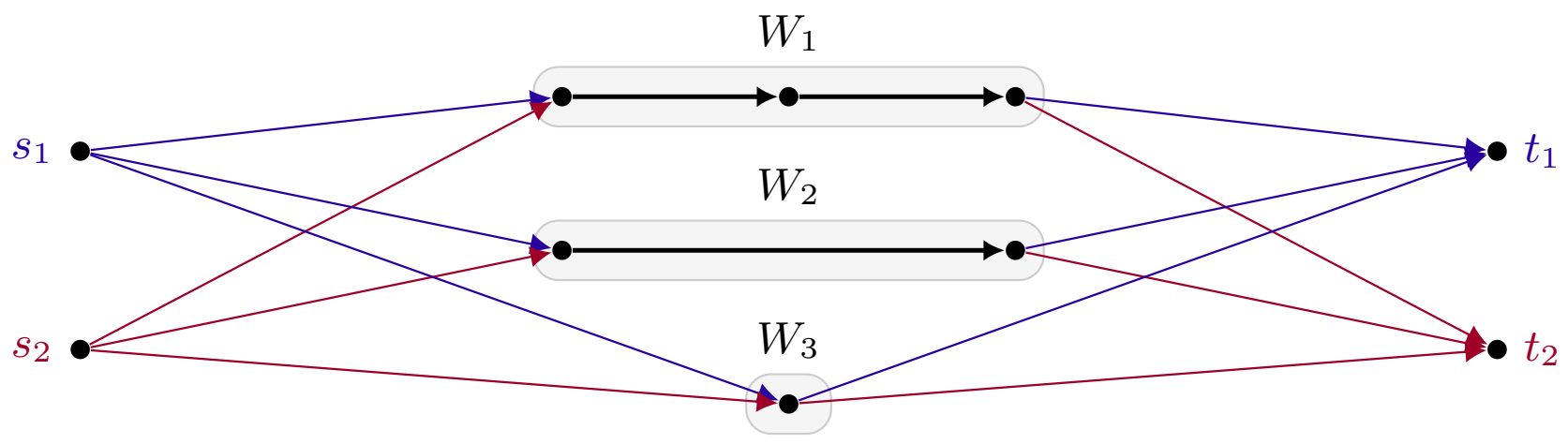

A \emph{parallel-channel topology} has distinct terminals (sources and sinks) $s_1,\ldots,s_k$, $t_1,\ldots,t_k$ and mutually vertex-disjoint directed
paths~$W_1,\ldots,W_\nu$ disjoint from the terminals, called \emph{channels}; a
channel may be a single vertex. Its only other edges are the
\emph{connectors}: one from each~$s_i$ to the first vertex of each~$W_q$,
and one from the last vertex of each~$W_q$ to each~$t_i$. The commodities
are~$(s_i,t_i)$ for $i\in[k]$; see~\cref{hier:fig-parallel-channels}.

\begin{lemma}\label{hier:parallel}
Every parallel-channel topology belongs to~$\mathcal S_3$.
\end{lemma}
\begin{proof}
For given capacities and demands, let~$a_{iq}$ be the minimum of the
capacities of the two connectors joining commodity~$i$ to channel~$q$.
Let~$c_q$ be the minimum capacity of an internal edge of~$W_q$, taking
$c_q=+\infty$ if~$W_q$ is a
single vertex. 
Every commodity path uses exactly one channel, so feasibility is
equivalent to the existence of $y_{iq}\ge0$ with
\begin{equation}\label{hier:allocation}
 \sum_q y_{iq}=d_i,\qquad y_{iq}\le a_{iq},\qquad
 \sum_i y_{iq}\le c_q.
\end{equation}
This is equivalent to a flow problem in an auxiliary network with edges from a source to each
commodity node~$i$ of capacity~$d_i$, from~$i$ to each channel node~$q$ of
capacity~$a_{iq}$, and from~$q$ to a sink of capacity~$c_q$. A cut with
commodity set~$S$ on the source side thus has capacity at least
$d([k]\setminus S)+\sum_q\min\{c_q,\sum_{i\in S}a_{iq}\}$, with equality
when each channel is placed optimally. By the max-flow/min-cut theorem,
\eqref{hier:allocation} is feasible if and only if
\begin{equation}\label{hier:channel-cuts}
 d(S)\le\sum_q\min\Bigl\{c_q,\sum_{i\in S}a_{iq}\Bigr\}
 \qquad\text{for }S\subseteq[k].
\end{equation}
Assume~\eqref{C3}, and fix~$S$. For each channel~$q$, take either an edge
of~$W_q$ of capacity~$c_q$ or, for every $i\in S$, a connector of
capacity~$a_{iq}$, whichever is cheaper; for a single-vertex channel only
the second option exists. The union~$F_S$ blocks every commodity in~$S$
and has capacity equal to the right-hand side
of~\eqref{hier:channel-cuts}, so~\eqref{C3} implies~\eqref{hier:channel-cuts}.
\qed
\end{proof}

\begin{corollary}\label{hier:ex2}
The second topology of~\cref{hier:fig-examples} lies
in~$\mathcal S_3\setminus\mathcal S_2$.
\end{corollary}
\begin{proof}
It is a parallel-channel topology with channels~$a$ and $(b,c)$, so it
lies in~$\mathcal S_3$. Let $e=(s_1,a)$, $f=(a,t_2)$, and $g=(b,c)$.
Blocking commodity~$1$ requires~$e$ and~$g$, and blocking commodity~$2$
requires~$f$ and~$g$. No cut can contain both~$e$ and~$f$. So a cut blocks at most one commodity, and then has
capacity at least $2=d_i$; thus~\eqref{C2} holds. The set $\{e,f,g\}$
blocks both commodities with capacity $3<4$, violating~\eqref{C3}.
\qed
\end{proof}

The violating set $\{e,f,g\}$ is an instance of the sets~$F_S$ in the
proof of~\cref{hier:parallel}, and it is not a cut. This is why
parallel channels need~\eqref{C3} rather than~\eqref{C2}.

\subsection{Vertex sums and supply-graph families}
Call a supply digraph \emph{\textup{(C$r$)}-sufficient} if
$(G,H)\in\mathcal S_r$ for every demand graph~$H$ on its vertices.
A \emph{vertex sum} of two digraphs~$G_1,G_2$ identifies one vertex
from each and leaves their edge sets disjoint.

\begin{lemma}\label{hier:vertex-sum}
For $r\in\{1,2,3,4\}$, vertex sums of~\textup{(C$r$)}-sufficient
supply digraphs are~\textup{(C$r$)}-sufficient; for $r=5$ this holds
trivially.
\end{lemma}
\begin{proof}
Let~$G$ be the vertex sum of~$G_1,G_2$ at~$v$, and suppose that~$(u,d)$ satisfies~\textup{(C$r$)}.  Discard zero-demand commodities; every remaining
commodity has a path, since~\eqref{C1} holds.
A simple path whose endpoints lie in one of~$G_1,G_2$ stays
in it, while a path between~$G_1-v$ and~$G_2-v$ passes through~$v$
exactly once.  Split each such commodity at~$v$, retaining its demand in
each piece; assign the others to the piece containing their endpoints,
where a commodity with an endpoint at~$v$ goes to the piece containing
its other endpoint. This gives
instances on~$G_1$ and~$G_2$; we show that the one on~$G_1$
satisfies~\textup{(C$r$)} (symmetric arguments hold for the instance on~$G_2$).

An edge set~$F\subseteq E(G_1)$ blocks a local commodity exactly when it blocks the corresponding original commodity; commodities having no local
segment remain unblocked.  Thus~\eqref{C3} is inherited.
For~\eqref{C4}, extend a binary length function by zero on~$G_2$.
The distance of each original commodity equals its local distance,
or zero if it has no segment in~$G_1$, so the local inequality is
again inherited.  For~\eqref{C1} and~\eqref{C2}, extend a vertex set~$X\subseteq V(G_1)$ by all of~$V(G_2)$ if $v\in X$, and leave it
unchanged otherwise.  Its outgoing edges are exactly~$\delta_{G_1}^+(X)$; the separated commodities for~\eqref{C1}, and
the blocked commodities for~\eqref{C2}, correspond to their local
segments.

By the assumed sufficiency, both instances are feasible.
At~$v$, join the two segment
flows of each split commodity, splitting path flows to match amounts. The
result is a feasible multiflow in~$G$, since~$G_1$ and~$G_2$ have disjoint edge sets.
\qed
\end{proof}

An \emph{edge-series-parallel} (ESP) digraph is obtained from a
single directed edge by series and parallel composition, with the
usual consistent terminal identifications.  Such digraphs are~\eqref{C2}-sufficient by~\cite[Thm.~3]{majthoub2026integer}.
A \emph{block} of an undirected multigraph is a maximal connected subgraph with no cut vertex, where a cut vertex is one whose removal (with its incident edges) disconnects the graph.
A \emph{block-ESP digraph} is a digraph in which every block of the
underlying undirected multigraph induces an ESP digraph for a suitable choice of terminals.
A \emph{cactus digraph} is an orientation
of an undirected cactus: each block is either an edge or a simple cycle.

\begin{corollary}\label{hier:cactus}
Cactus digraphs and block-ESP digraphs are~\eqref{C2}-sufficient.
\end{corollary}
\begin{proof}
By~\cref{hier:vertex-sum}, it suffices to consider individual blocks;
disconnected components can be treated separately, since~\eqref{C1}
excludes positive demand between them.  The assertion for ESP blocks
is the cited theorem.  In an orientation of a simple undirected cycle,
two distinct directed paths with the same endpoints, if they exist,
cover the cycle.  The block is then the parallel composition of these
two paths, hence an ESP digraph.  Otherwise every ordered pair has
at most one simple directed path and~\cref{hier:unique} applies.
\qed
\end{proof}

Neither family contains the other: three internally disjoint directed
paths of length at least two with common endpoints form an ESP block whose underlying graph
is not a cactus, whereas a consistently directed cycle is a cactus
block and is not ESP.

\section{Complexity of checking and recognition}
\label{hier:appendix-complexity}

\subsection{Checking the conditions}
% \label{hier:appendix-complexity}
Here the input is a weighted instance with rational capacities and
demands encoded in binary. Recognizing sufficient topologies is treated
in~\cref{hier:appendix-recognition}.

\begin{proposition}\label{hier:easy-checks}
Condition~\eqref{C1} can be checked with at most~$4^k$ minimum-cut
computations.  Condition~\eqref{C5} can be checked in polynomial time,
without a bound on the number of commodities.
\end{proposition}
\begin{proof}
There are at most~$2k$ terminals. For each set~$T$ of terminals, all
sets~$X$ containing exactly the terminals in~$T$ have the same right-hand
side in~\eqref{C1}. If it is positive, one minimum-cut computation, with
the terminals in~$T$ merged into a source and the others into a sink,
finds such an~$X$ of minimum capacity. This also gives the
fixed-parameter tractability observed
in~\cite[Cor.~3]{Poremba_Shepherd_2025}. For~\eqref{C5}, solve the
edge-flow feasibility linear program, which has $k|E|$ variables.
\qed
\end{proof}

\begin{theorem}\label{hier:hard-checks}
For each $r\in\{2,3,4\}$, deciding whether a two-commodity instance
satisfies~\textup{(C$r$)} is strongly coNP-complete.
\end{theorem}
\begin{proof}
\emph{Membership in coNP.} A violating vertex set, edge set, or binary length
function is a polynomial-size certificate for the complement, since blocked
commodities and binary distances are computable in polynomial time.

\emph{Hardness for~\eqref{C3}.} Bentz~\cite{Bentz_11} proved that it is
NP-hard to decide whether an unweighted directed acyclic graph has an
edge set of size at most~$B$ blocking two given terminal pairs.
Let~$\lambda_i$ be the minimum size of an $s_i$--$t_i$ cut. The answer
is no if $B<\max\{\lambda_1,\lambda_2\}$, and yes if
$B\ge\lambda_1+\lambda_2$, as the union of two minimum cuts shows. Since
the~$\lambda_i$ are computable in polynomial time, the problem remains
NP-hard on instances with
$\max\{\lambda_1,\lambda_2\}\le B\le\lambda_1+\lambda_2-1$, and we reduce
from these. Give every edge capacity~$1$, and set $d_1:=\lambda_1$ and
$d_2:=B+1-\lambda_1$, so that $1\le d_i\le\lambda_i$ and
$d_1+d_2=B+1$. An edge set blocking only commodity~$i$ has capacity at
least $\lambda_i\ge d_i$. Hence~\eqref{C3} fails exactly when an edge
set of size at most~$B$ blocks both commodities. All numbers are bounded
by~$|E|$, so the hardness is strong.

\emph{From checking~\eqref{C3} to checking~\eqref{C2}.}
Let $d([k])=\sum_i d_i$ and replace each edge~$e=(v,w)$ by~$v\to z_e\to w$, with capacities~$d([k])$ and~$u_e$, respectively.
A set~$F\subseteq E$ corresponds to the cut induced by
$Z_F:=\{z_e:e\in F\}$, which consists of the second
halves of the edges in~$F$.  This cut has capacity~$u(F)$ and blocks
the same commodities as~$F$.  Conversely, a cut violating~\eqref{C2}
cannot contain a first half, since that alone has capacity~$d([k])$.
Its outgoing second halves project to an edge set of the original
graph of the same capacity blocking the same commodities.
Thus, the subdivided instance satisfies~\eqref{C2} if and only if
the original satisfies~\eqref{C3}.

\emph{From checking~\eqref{C3} to checking~\eqref{C4}.}
Add, for each commodity, a new
edge $\hat e_i=(s_i,t_i)$ of capacity~$0$. For a binary length vector~$\ell$
on the augmented graph, let $F:=\{e\in E(G):\ell_e=1\}$. The new edges
make every distance at most~$1$, and a commodity not blocked by~$F$ has
distance~$0$. Hence
\[
 \sum_eu_e\ell_e=u(F),\qquad
 \sum_id_i\dist_\ell(s_i,t_i)\le d(K(F)),
\]
with equality when every~$\hat e_i$ has length~$1$. Thus~\eqref{C4} for the
augmented instance is equivalent to~\eqref{C3} for the original.

Both reductions keep two commodities and polynomially bounded numbers.
\qed
\end{proof}

\subsection{Recognizing sufficient topologies}
\label{hier:appendix-recognition}

We prove the two parts of~\cref{hier:recognition} separately.
% \paragraph{Membership in coNP.}
\begin{lemma}[Membership]\label{lem:membership}
    For each $r\in\{1,2,3,4\}$, deciding whether a topology belongs
to~$\mathcal S_r$ is in coNP, for any number of commodities.
\end{lemma}
\begin{proof}
A weighting that satisfies~\textup{(C$r$)} but is infeasible would
certify nonmembership. However, checking~\textup{(C$r$)} is coNP-hard
for $r\ge2$ (\cref{hier:hard-checks}), and~\eqref{C1} is only known to
be checkable efficiently for fixed~$k$. We therefore certify
nonmembership by a single inequality.

\emph{Setup.} A commodity without a path has demand zero
under~\eqref{C1}, and hence under each of the four conditions, so
deleting it does not affect membership in~$\mathcal S_r$. Assume that
every commodity has a path, and let $N:=|E|+k\ge1$. The feasible
capacity--demand vectors form the cone $\mathcal F\subseteq\R^N$
generated by
\[
 \{(\ind{\{e\}},0):e\in E\}\ \cup\
 \{(\ind{E(P)},\ind{\{i\}}):i\in[k],\ P\in\mathcal P_i\},
\]

where~$\ind{\{e\}}$ and~$\ind{\{i\}}$ are unit vectors of the capacity
and demand coordinates, respectively.  The first generators represent
unused capacity and the second unit path flows.  These $0/1$ vectors
span~$\R^N$ (take all~$(\ind{\{e\}},0)$ and one path vector per
commodity), so~$\mathcal F$ is full-dimensional. The inequality
$\sum_eu_e\ell_e\ge\sum_iq_id_i$ is valid for~$\mathcal F$ if and only if
$\ell\ge0$ and $q_i\le\dist_\ell(s_i,t_i)$ for all~$i$.

\emph{Facets.} Let~$\mathcal Q_r$ be the cone defined by nonnegativity
and the inequalities of~\textup{(C$r$)}. These are valid for~$\mathcal F$,
so $\mathcal F\subseteq\mathcal Q_r$, and $(G,H)\in\mathcal S_r$ if and
only if $\mathcal F=\mathcal Q_r$. This holds if and only if every facet
inequality of~$\mathcal F$ is a positive multiple of a defining
inequality of~$\mathcal Q_r$. If it is, $\mathcal Q_r$ lies in every facet
halfspace, hence in~$\mathcal F$. Conversely, if
$\mathcal F=\mathcal Q_r$, Farkas' lemma writes each facet normal as a
nonnegative combination of defining normals of~$\mathcal Q_r$. These lie
in the dual cone~$\mathcal F^*$, and a facet normal spans an extreme ray
of~$\mathcal F^*$, so each normal with a positive coefficient is a
positive multiple of it.

\emph{Certificate.} Nonmembership is thus certified by a facet inequality
of~$\mathcal F$ that \emph{matches} no defining inequality
of~$\mathcal Q_r$, that is, is a positive multiple of none. Such a certificate consists of a nonzero integer normal~$(\ell,-q)$
and~$N-1$ linearly independent generators of~$\mathcal F$ that are tight
for the inequality
$
 \sum_e u_e\ell_e\ge\sum_i q_i d_i.
$
Each generator can be encoded by an edge or a simple commodity path.
An integer facet normal of polynomial encoding length exists: take a
cofactor vector of the matrix of the~$N-1$ tight generators, choosing
its overall sign so that the inequality is valid.  Hadamard's inequality
bounds the encoding length of each
entry by~$O(N\log N)$ for $N\ge2$; the case~$N=1$ is immediate.
Validity of the inequality is equivalent to
$
 \ell\ge0,\qquad
 q_i\le\dist_\ell(s_i,t_i)\quad\text{for }i\in[k],
$
which is checked by shortest-path computations.  Tightness and rank are
also polynomially checkable, and together with validity certify that
the inequality defines a facet.

\emph{Matching.} Coordinate inequalities are matched directly. Each
inequality of \textup{(C$r$)} has capacity part~$\ind{F'}$ for some
$F'\subseteq E$, and is trivial if $F'=\emptyset$, since all commodities
have paths. So a match requires all positive entries of~$\ell$ to equal
some $L>0$; then $F'=F:=\{e:\ell_e=L\}$, with factor~$L$.
\begin{itemize}
\item For~\eqref{C3} and~\eqref{C4}, a match exists if and only if
  $q=L\,\ind{K(F)}$, respectively $q_i=L\dist_{\ind{F}}(s_i,t_i)$,
  for all~$i$.
\item For~\eqref{C2}, a match additionally requires $F=\delta^+_G(X)$ for
  some~$X$. Let~$R$ be the set reachable in~$G-F$ from the tails
  of~$F$. Such an~$X$ exists if and only if~$R$ contains no head of~$F$:
  any such~$X$ contains the tails, avoids the heads, and is closed under
  edges outside~$F$, so $R\subseteq X$; conversely, $X=R$ works.
\item For~\eqref{C1}, a match requires $q_i\in\{0,L\}$ and a set~$X$ with
  $\delta^+_G(X)=F$ and $\delta^+_H(X)=\{i:q_i=L\}$. Force the tails
  of~$F$ and the sources~$s_i$ with $q_i=L$ into~$X$, and the heads
  of~$F$ and the sinks~$t_i$ with $q_i=L$ out of it. The set~$X$ must
  contain the head of every edge of $E\setminus F$ whose tail it
  contains, and~$t_i$ whenever it contains~$s_i$ with $q_i=0$. Such
  an~$X$ exists if and only if the set reachable from the forced-in
  vertices along these edges and pairs avoids the forced-out vertices;
  that set is then one.
\end{itemize}
All these tests run in polynomial time.
\qed
\end{proof}

For hardness, we reduce from the \emph{useful edge problem}: given a
digraph~$D$, distinct vertices~$s,t$, and an edge~$e=(w,v)$, is there an
$s$--$t$ path that uses~$e$? This problem is
NP-complete~\cite{Fortune_Hopcroft_Wyllie_80}, and Poremba and
Shepherd~\cite[Sect.~5]{Poremba_Shepherd_2025} also reduce from it.

In the following proof we use the notion of \emph{$b$-transshipments (transshipments)}. Let $b\in\R^V$ with $\sum_{v\in V}b(v)=0$,
and write $b(X):=\sum_{v\in X}b(v)$. A \emph{$b$-transshipment} is a
flow~$f$ in~$G$ with $0\le f\le u$ whose net outflow at every
vertex~$v$ is~$b(v)$. It exists if and only if
\[
 b(X)\le u(\delta^+_G(X))\qquad\text{for every }X\subseteq V.
\]
This follows from the max-flow/min-cut theorem after adding a source
with an edge of capacity~$b(v)$ to every vertex with $b(v)>0$, and a
sink with an edge of capacity~$-b(v)$ from every vertex with $b(v)<0$.

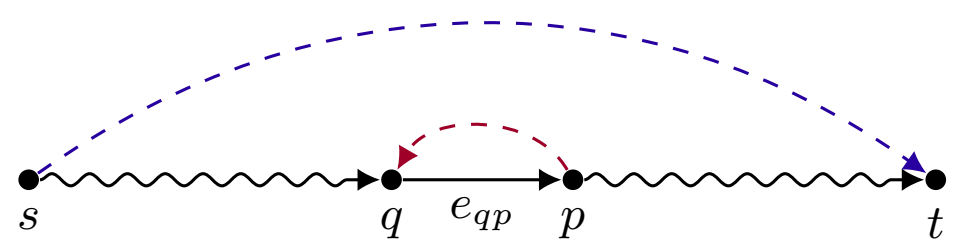
\begin{figure}[htb]
\centering
\begin{tikzpicture}[x=1cm,y=1cm]
  \node[solidnode,label=below:$s$] (s) at (0,0) {};
  \node[solidnode,label=below:$q$] (q) at (2.4,0) {};
  \node[solidnode,label=below:$p$] (p) at (3.6,0) {};
  \node[solidnode,label=below:$t$] (t) at (6,0) {};
  \draw[wavy] (s) -- (q);
  \draw (q) edge[dedge] node[below]{$e_{qp}$} (p);
  \draw[wavy] (p) -- (t);
  \draw[demandedge,ccolor1] (s) to[bend left=35] (t);
  \draw[demandedge,ccolor2] (p) to[bend right=60] (q);
\end{tikzpicture}
\caption{The proof of~\cref{hier:reverse-arc}. Wavy lines are paths
in~$G$; dashed edges are demands. After the edge $e_{qp}=(q,p)$ is added, an
$s$--$t$ path through~$e_{qp}$ consists of vertex-disjoint $s$--$q$ and
$p$--$t$ paths.}
\label{hier:fig-reverse-arc}
\end{figure}

\begin{lemma}\label{hier:reverse-arc}
Let~$(G,H)$ have the two commodities $(s,t)$ and~$(p,q)$. If~$G$ has no
vertex-disjoint paths $s$--$q$ and $p$--$t$, then
$(G,H)\in\mathcal S_1$.
\end{lemma}

\begin{proof}
Let $(u,d)$ satisfy~\eqref{C1}; we construct a feasible multiflow. Put
\[
b:=d_1(\ind{\{s\}}-\ind{\{t\}})+d_2(\ind{\{p\}}-\ind{\{q\}}).
\]
For every $X\subseteq V$,
$b(X)$ is at most the demand separated by~$X$, so
$b(X)\le d(\delta^+_H(X))\le u(\delta^+_G(X))$ by~\eqref{C1}. Hence
there is a $b$-transshipment~$f$. It moves $d_1+d_2$ units from~$s$
and~$p$ to~$t$ and~$q$, but it may send flow from~$s$ to~$q$ or from~$p$
to~$t$.

To separate the commodities, add the edge $e_{qp}:=(q,p)$ with $d_2$ units of
flow on it. The resulting flow has net outflow $d_1$ at~$s$, $-d_1$
at~$t$, and zero elsewhere, so it decomposes into $s$--$t$ paths of total
value~$d_1$ and cycles. Any simple $s$--$t$ path through~$e_{qp}$ would consist of an
$s$--$q$ path, the edge~$e_{qp}$, and a $p$--$t$ path, and these two paths would have to be
vertex-disjoint for the simplicity of the whole path
(\cref{hier:fig-reverse-arc}). By assumption there are no such vertex-disjoint paths, so~$e_{qp}$ does not lie in any $s$--$t$ path; all $d_2$ units on~$e_{qp}$ lie on cycles. Removing~$e_{qp}$ from each of these
cycles leaves $p$--$q$ paths of total value~$d_2$. Together with the
$s$--$t$ paths, they have flow at most $f\le u$ on every edge, so they form a
feasible multiflow.
\qed
\end{proof}

\begin{figure}[!b]
\centering
\begin{tikzpicture}[x=1cm,y=1cm]
  \draw[rounded corners=5pt,fill=black!4,draw=black!20]
    (-.4,-.55) rectangle (7,.55);
  \node[font=\small] at (.2,.8) {$D-e$};
  \node[solidnode,label=left:$s$] (s) at (0,0) {};
  \node[solidnode,label=below:$w$] (w) at (2,0) {};
  \node[solidnode,label=below:$v$] (v) at (4.6,0) {};
  \node[solidnode,label=right:$t$] (t) at (6.6,0) {};
  \draw[wavy,black!50] (s) -- node[above=1pt]{\footnotesize$P_{sw}$} (w);
  \draw[wavy,black!50] (v) -- node[above=1pt]{\footnotesize$P_{vt}$} (t);
  \draw[dedge,black!25] (w) -- node[below]{\footnotesize$e$} (v);
  % gadget: the fourth topology of Fig. 2, a,z,b,c = s_1,t_1,s_2,t_2
  \node[solidnode,label=left:$a$]  (a) at (2.5,1.1) {};
  \node[solidnode,label=right:$z$] (z) at (4.1,1.1) {};
  \node[solidnode,label=left:$b$]  (b) at (2.5,2.5) {};
  \node[solidnode,label=right:$c$] (c) at (4.1,2.5) {};
  \draw (a) edge[dedge,dotted] (b)
        (c) edge[dedge,dotted] (z)
        (z) edge[dedge,dotted] (a)
        (a) edge[dedge] (c)
        (b) edge[dedge] (z)
        (b) edge[dedge] (c);
  \draw[dedge] (w) -- (a);
  \draw[dedge] (z) -- (v);
  \draw[demandedge,ccolor1] (s) to[bend right=25] (t);
  \draw[demandedge,ccolor2] (b) to[bend left=40] (c);
\end{tikzpicture}
\caption{The topology of the hardness proof. The edge~$e$ (light) is
deleted, and the gadget on $a,b,c,z$ is attached through $(w,a)$
and~$(z,v)$. The gadget is the fourth topology
of~\cref{hier:fig-examples}, with $a,z,b,c$ in the roles of
$s_1,t_1,s_2,t_2$ and the same solid and dotted edges. Dashed edges are
the demands. The paths $P_{sw}$ and~$P_{vt}$ exist only in the second
case of the proof.}
\label{hier:fig-hardness}
\end{figure}
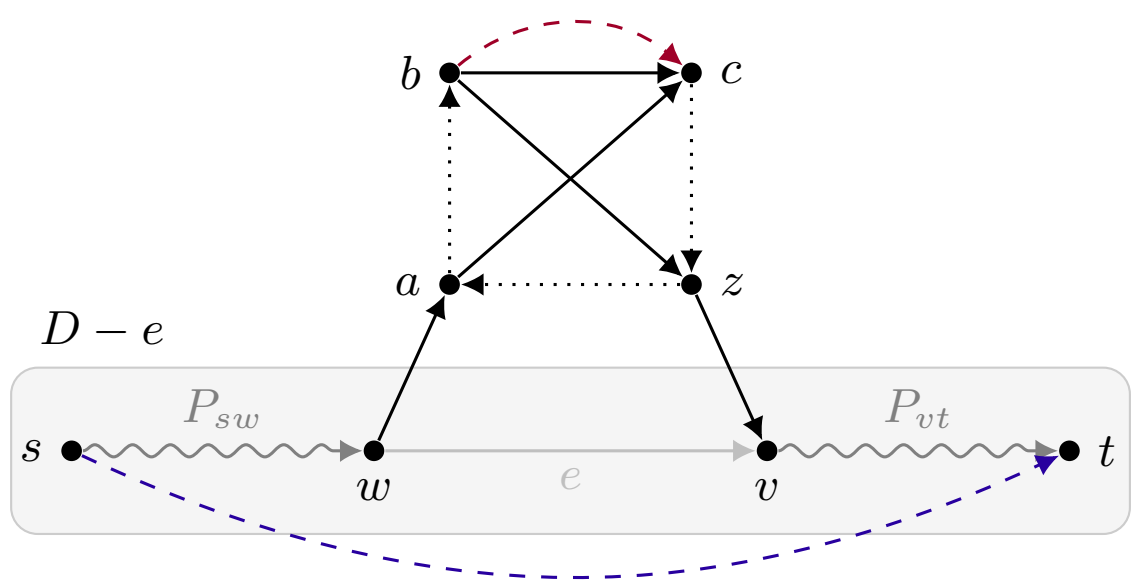

\begin{lemma}\label{lem:hardness}
For each $r\in\{1,2,3,4\}$, deciding whether a topology belongs
to~$\mathcal S_r$ is coNP-hard already for two-commodity topologies with two-matching demands.   
\end{lemma}
\begin{proof}
Let~$(D,s,t,e)$ with $e=(w,v)$ be an instance of the useful edge problem. Build~$G$ from~$D$ by deleting~$e$ and attaching the gadget
of~\cref{hier:fig-c4} through the edges $(w,a)$ and $(z,v)$ (see~\cref{hier:fig-hardness}). The
demand graph~$H$ has the commodities $(s,t)$ and $(b,c)$, so the four
terminals are distinct. Every path that enters the gadget uses $(w,a)$,
and every path that leaves it uses $(z,v)$.

We show that $(G,H)\in\mathcal S_1$ if no $s$--$t$ path in~$D$ uses~$e$,
and $(G,H)\notin\mathcal S_4$ otherwise. This proves hardness for every
class between~$\mathcal S_1$ and~$\mathcal S_4$, in particular for
each~$\mathcal S_r$ with $r\le4$.

\emph{No $s$--$t$ path in~$D$ uses~$e$.} By~\cref{hier:reverse-arc}, it suffices to show that~$G$ has no vertex-disjoint paths~$Q_1$ from~$s$ to~$c$ and~$Q_2$ from~$b$ to~$t$. Suppose such paths exist. Since $Q_1$ enters the gadget only once, via $(w,a)$, its prefix is an $s$--$w$ path in $D-e$. Similarly, $Q_2$ exits the gadget only once, via $(z,v)$, and its suffix is a $v$--$t$ path in $D-e$. These two paths are vertex-disjoint, so together with~$e$ they yield an $s$--$t$ path in~$D$ using~$e$, a contradiction.

\emph{An $s$--$t$ path $P=P_{sw}\cdot e\cdot P_{vt}$ in~$D$ uses~$e$.}
Let $d=(1,2)$. Assign capacities as in~\cref{hier:fig-c4}: set capacity~$1$ on solid edges $(a,c)$, $(b,z)$, $(b,c)$ and capacity~$3$ on dotted edges. Also set capacity~$3$ on $(w,a)$, $(z,v)$, and all edges of $P_{sw}$ and~$P_{vt}$; set capacity~$0$ on all remaining edges.

First, we observe that \eqref{C4} holds for this weighting. Let $\ell\in\{0,1\}^{E(G)}$, let $\kappa$ be the total
$\ell$-length on capacity-$3$ edges, and let $\ell_{ac}, \ell_{bz}, \ell_{bc}$ be its
values on $(a,c)$, $(b,z)$, $(b,c)$. Commodity~$1$ may use $P_{sw}$,
$(w,a)$, one of $a,c,z$, $a,b,z$, or $a,b,c,z$, then
$(z,v)$ and $P_{vt}$. Commodity~$2$ may use $(b,c)$ or 
$b,z,a,c$. Each route uses any capacity-$3$ edge at most
once, so with $\ell\ge0$,
\[
\dist_\ell(s,t)\le \kappa + \min\{\ell_{ac}, \ell_{bz}, \ell_{bc}\},\qquad
\dist_\ell(b,c)\le \kappa + \min\{\ell_{bc},\ell_{ac}+\ell_{bz} \}.
\]
By \eqref{hier:eq-gadget},
\[
\dist_\ell(s,t)+2\dist_\ell(b,c)\le 3 \kappa+\ell_{ac}+ \ell_{bz}+\ell_{bc}=\sum_e u_e\ell_e.
\]

Now for infeasibility, note that positive capacity lies only on~$P_{sw}$,
$P_{vt}$, $(w,a)$, $(z,v)$ and the gadget edges. Using only these
edges, every $s$--$t$ path crosses the gadget from~$a$ to~$z$, and
every $b$--$c$ path is the edge $(b,c)$ or the path
$b,z,a,c$. Give length~$1$ to $(a,c)$ and $(b,z)$, length~$2$ to
$(b,c)$, and length~$0$ to all other edges. Then every such $s$--$t$ path
has length at least~$1$, and every such $b$--$c$ path has length at
least~$2$. A feasible multiflow would therefore need
$\sum_eu_e\ell_e\ge1\cdot1+2\cdot2=5$, but $\sum_eu_e\ell_e=4$.
\qed
\end{proof}

Now we are in the position to prove \cref{hier:recognition}.

\begin{proof}[of~\cref{hier:recognition}]
    This follows by combining \cref{lem:hardness,lem:membership}. 
\end{proof}

\section{A weight-preserving reduction for blocking-set violations}
\label{char:blocking-reduction}

The contraction in the proof of~\cref{char:main} may not preserve infeasibility of the
given weights. In the instance of~\cref{hier:ex4}, the
unique minimum-length integral witness is 
$(\ell_e,\ell_f,\ell_g)=(1,1,2)$ (and $0$ on the dotted edges).
 By the reduction in~\cref{hier:appendix-families},
any minimum witness is $0$ on dotted edges and has
total length at most~$4$; if any of $\ell_e,\ell_f,\ell_g$ is $0$, then the
distance formulas in~\cref{hier:ex4} imply~\eqref{C5}. Among
 positive triples with sum at most~$4$, only $(1,1,2)$
violates~\eqref{C5}.

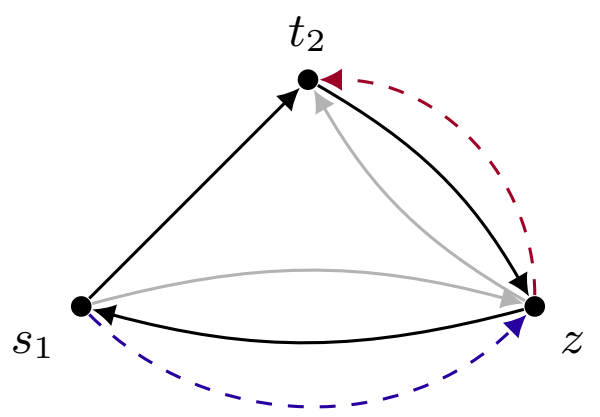
\begin{figure}[tb]
\centering
\begin{tikzpicture}[x=1cm,y=1cm]
  \node[solidnode,label=below left:$s_1$] (s) at (0,0) {};
  \node[solidnode,label=above:$t_2$] (t) at (1.5,1.5) {};
  \node[solidnode,label=below right:$z$] (z) at (3,0) {};
  \draw (s) edge[dedge] (t)
        (t) edge[dedge,bend left=15] (z)
        (z) edge[dedge,bend left=15] (s)
        (s) edge[dedge,black!30,bend left=15] (z)
        (z) edge[dedge,black!30,bend left=15] (t);
  \draw[demandedge,ccolor1] (s) to[bend right=45] (z);
  \draw[demandedge,ccolor2] (z) to[bend right=45] (t);
\end{tikzpicture}
\caption{The instance of~\cref{hier:ex4} after contracting
$f=(s_2,t_1)$ into~$z$. Solid edges are supply edges; dashed edges are
demands. With the original weights the instance is feasible, but
deleting the two gray edges leaves the bad triangle.}
\label{char:fig-contracted}
\end{figure}

For this witness, $m_1=1$ and $m_2=2$. The only short cross-commodity
path is $f=(s_2,t_1)$ with length $1<m_2$, so the proof
contracts $f$ to a vertex~$z$ (\cref{char:fig-contracted}). With the
original weights, the contracted instance becomes feasible: commodity~1 
uses $(s_1,z)$, and commodity~2 sends one unit on $(z,t_2)$
and one on $z,s_1,t_2$. Deleting $(s_1,z)$ and
$(z,t_2)$ still leaves the bad triangle as a relevant minor. However, when~\eqref{C3}
fails, a different contraction keeps the original weights infeasible.

\begin{lemma}[Cross-reachability of blocked commodities]
\label{char:cross-blocked}
Let a multiflow instance satisfy~\eqref{C1}, and let $F\subseteq E(G)$
violate~\eqref{C3}. Then there are distinct positive-demand commodities
$i,j\in K(F)$ such that $G-F$ contains an $s_i$--$t_j$ path.
\end{lemma}
\begin{proof}
Let $J:=\{i\in K(F):d_i>0\}$ and $X:=\bigcup_{i\in J}X_i(F)$.
Then, by~\eqref{eq:canonical-cut}, $\delta^+_G(X)\subseteq F$, and~$X$ contains
every~$s_i$ with $i\in J$. If~$X$ contained no~$t_j$ with $j\in J$,
then~\eqref{C1} would give
$u(F)\ge u(\delta^+_G(X))\ge d(J)=d(K(F))$, a contradiction. As a consequence,
$t_j\in X_i(F)$ for some $i,j\in J$, and $i\ne j$ because
$t_j\notin X_j(F)$.
\qed
\end{proof}

\begin{proposition}\label{char:c3-reduction}
Let~$(G,H,u,d)$ be a two-matching instance that satisfies~\eqref{C1} and
violates~\eqref{C3}. After possibly swapping the commodities, there is an
$s_1$--$t_2$ path whose contraction gives a two-path instance, with the
original weights on the surviving edges, that satisfies~\eqref{C1} and
violates~\eqref{C3}.
\end{proposition}
\begin{proof}
Choose $F\subseteq E(G)$ violating~\eqref{C3}.
By~\cref{char:cross-blocked}, both commodities are blocked by~$F$. After swapping them if necessary, $G-F$ contains an $s_1$--$t_2$ path~$P$. The
path~$P$ avoids~$t_1$, since commodity~$1$ is blocked by~$F$, and it
avoids~$s_2$, since otherwise its part from~$s_2$ to~$t_2$ would
avoid~$F$. So~\cref{char:safe} applies: contracting~$P$ into a new
vertex~$z$ gives a relevant minor with the two-path demands $(s_2,z)$
and $(z,t_1)$.

Keep the original capacities and demands on the surviving edges, and
let~$\bar F$ be the set of edges of~$F$ that do not become loops. Then
$u(\bar F)\le u(F)<d_1+d_2$. We show that~$\bar F$ blocks both
commodities in the minor. A $z$--$t_1$ path avoiding~$\bar F$ comes from
a path in~$G$ that starts at a vertex~$p$ of~$P$, ends at~$t_1$,
avoids~$F$, and meets~$P$ only at~$p$. Together with~$P[s_1,p]$, it
would give an $s_1$--$t_1$ path avoiding~$F$. Similarly, an $s_2$--$z$
path avoiding~$\bar F$, followed by a part of~$P$, would give an
$s_2$--$t_2$ path avoiding~$F$. Both contradict $1,2\in K(F)$, so~$\bar F$
violates~\eqref{C3} in the minor.

Finally, a vertex set~$Y'$ of the minor corresponds to the vertex
set~$Y$ of~$G$ obtained by replacing~$z$ with~$V(P)$. The edges
leaving~$Y'$ are exactly the edges leaving~$Y$, and~$Y'$ and~$Y$
separate the same commodities. So the minor satisfies~\eqref{C1}.
\qed
\end{proof}

Combined with~\cref{char:ps}, this proves~\cref{char:main} for deceptive
weights that violate~\eqref{C3}, without the witness argument
of~\cref{char:sec-cross,char:sec-layer}. That argument is needed only for
deceptive weights that satisfy~\eqref{C3}.

\section{Consequences of~\cref{char:main}}
Define the \emph{fractional flow-cut gap}~$\beta_f(G,H)$ as the infimum of all
$\beta\ge1$ such that any choice of nonnegative finite weights
satisfying~\eqref{C1} has a routing under capacities~$\beta u$.  Define
the \emph{integral flow-cut gap}~$\beta_{\mathrm{int}}(G,H)$ in the same
way, requiring an integral routing and restricting the original
capacities and demands to integers.

\begin{corollary}\label{char:gap}
For every directed two-commodity topology,
\[
 \beta_f(G,H)=\beta_{\mathrm{int}}(G,H)=
 \begin{cases}
 1,&\text{if~$(G,H)$ is cut-sufficient},\\
 2,&\text{otherwise}.
 \end{cases}
\]
\end{corollary}
\begin{proof}
The cut condition allows each commodity to be routed individually
within capacity~$u$.  Superposing these routings gives a multiflow
within capacity~$2u$, integral when~$u,d$ are integral.
If the topology is cut-sufficient, the Nagamochi--Ibaraki
theorem~\cite{Nagamochi_Ibaraki_89}, or \cref{int:main} (since $\mathcal S_1\subseteq\mathcal S_4$) gives both gaps equal to one.
Otherwise~\cref{char:main} gives a forbidden relevant minor.
Its unit capacities and demands require scaling by at least two,
because both commodity paths use a common edge.  Extend these weights
to the original topology.  The extension is integral and satisfies~\eqref{C1}, while any routing under capacities~$\beta u$ would project
to such a routing on the minor.  No $\beta<2$ is therefore possible.
\qed
\end{proof}

\begin{remark}\label{char:unit-certificate}
If a directed two-commodity topology is not cut-sufficient, it admits
deceptive weights with demands in~$\{0,1\}$ and capacities
in~$\{0,1,2\}$: as in the proof of~\cref{char:gap}, extend the unit
weights of a forbidden relevant minor (\cref{char:extension}, with
$D=2$). This confirms the unit-demand certificate property conjectured
in~\cite[Sect.~6]{Poremba_Shepherd_2025}.
\end{remark}

\begin{remark}\label{char:bridge}
Each forbidden topology lies in~$\mathcal S_2\setminus\mathcal S_1$
by~\cref{hier:unique}. Hence a two-commodity topology is
cut-sufficient if and only if it has no relevant minor in
$\mathcal S_2\setminus\mathcal S_1$. This concerns topologies:
certifying infeasibility of particular weights may require~\eqref{C5},
and the contraction above need not preserve those weights'
infeasibility.
\Cref{char:blocking-reduction} gives an example, together with a contraction that does preserve it when~\eqref{C3} fails.
\end{remark}

\section{Gaps at the higher levels}
\label{gap:appendix}

By~\cref{char:gap}, the cut condition has gap one or two on every
two-commodity topology. Here we show that the higher levels allow every
gap value $2-1/r$
(\cref{gap:sec:A,gap:sec:B,gap:sec:consequences}).

\subsection{Gaps at the higher levels}

For weights $(u,d)$ satisfying~\eqref{C1}, let $\gamma(u,d)$ be the least
factor~$\gamma$ such that $(G,H,\gamma u,d)$ is feasible. For
$\lambda\in\{1,2,3,4\}$, the \emph{fractional gap at level~$\lambda$},
denoted~$\beta_\lambda(G,H)$, is the larger of~$1$ and the supremum of $\gamma(u,d)$
over all weights satisfying~\textup{(C$\lambda$)}. The \emph{integral gap}
$\beta_\lambda^{\mathrm{int}}(G,H)$ is defined in the same way, with integral
weights and integral routings. Thus $\beta_1=\beta_f$ and
$\beta_1^{\mathrm{int}}=\beta_{\mathrm{int}}$ in~\cref{char:gap}.
Stronger conditions allow fewer weights, and routing the two commodities
separately, as in the proof of~\cref{char:gap}, costs at most a factor
two. Hence, for two commodities,
\begin{equation}\label{gap:upper}
 1\le\beta_4\le\beta_3\le\beta_2\le\beta_1\le2
 \quad\text{and}\quad
 1\le\beta_4^{\mathrm{int}}\le\beta_3^{\mathrm{int}}
 \le\beta_2^{\mathrm{int}}\le\beta_1^{\mathrm{int}}\le2.
\end{equation}
By~\cref{int:main}, a \textup{(C$\lambda$)}-sufficient topology has both gaps
equal to one. The table summarizes the results for the two families
below, where $r\ge2$. The entry $[2-1/r,2]$ gives proved bounds, not an
exact value.
\begin{center}
\begin{tabular}{lccc}
\hline
 & Level 2 & Level 3 & Level 4\\
\hline
$\beta_\lambda(A_r)$ & $2-1/r$ & $2-1/r$ & $1$\\
$\beta_\lambda(B_r)$ & $[2-1/r,\,2]$ & $2-1/r$ & $2-1/r$\\
$\beta_\lambda^{\mathrm{int}}(A_r)$ & $2$ & $2$ & $1$\\
$\beta_\lambda^{\mathrm{int}}(B_r)$ & $2$ & $2$ & $2$\\
\hline
\end{tabular}
\end{center}

\subsection{The forced-chain family}\label{gap:sec:A}

\begin{figure}[!b]
\centering
\begin{tikzpicture}[x=1cm,y=1cm]
  \foreach \i in {1,2,3}{
    \node[solidnode,label=below:$a_{\i}$] (a\i) at ({2.4*(\i-1)},0) {};
    \node[solidnode,label=below:$b_{\i}$] (b\i) at ({2.4*(\i-1)+1},0) {};
    \draw (a\i) edge[dedge] node[above]{$e_{\i}$} (b\i);
  }
  \node[solidnode,label=left:$s_1$]  (s1) at (-1.2,0) {};
  \node[solidnode,label=right:$t_1$] (t1) at (7,0) {};
  \node[solidnode,label=left:$s_2$]  (s2) at (-1.2,1.8) {};
  \node[solidnode,label=right:$t_2$] (t2) at (7,1.8) {};
  \draw (s1) edge[dedge] (a1) (b1) edge[dedge] (a2)
        (b2) edge[dedge] (a3) (b3) edge[dedge] (t1);
  \foreach \i in {1,2,3}{
    \draw (s2) edge[dedge] (a\i) (b\i) edge[dedge] (t2);
  }
  \draw[demandedge,ccolor1] (s1) to[bend right=35] (t1);
  \draw[demandedge,ccolor2] (s2) to[bend left=25] (t2);
\end{tikzpicture}
\caption{The topology~$A_3$. Solid edges are supply edges; dashed edges
are demands. Crossings without a vertex are not junctions. In the
lower-bound weighting, the low edges~$e_j$ have capacity~$1$, all other
edges capacity~$3$, and $(d_1,d_2)=(1,2)$.}
\label{gap:fig:Ar}
\end{figure}
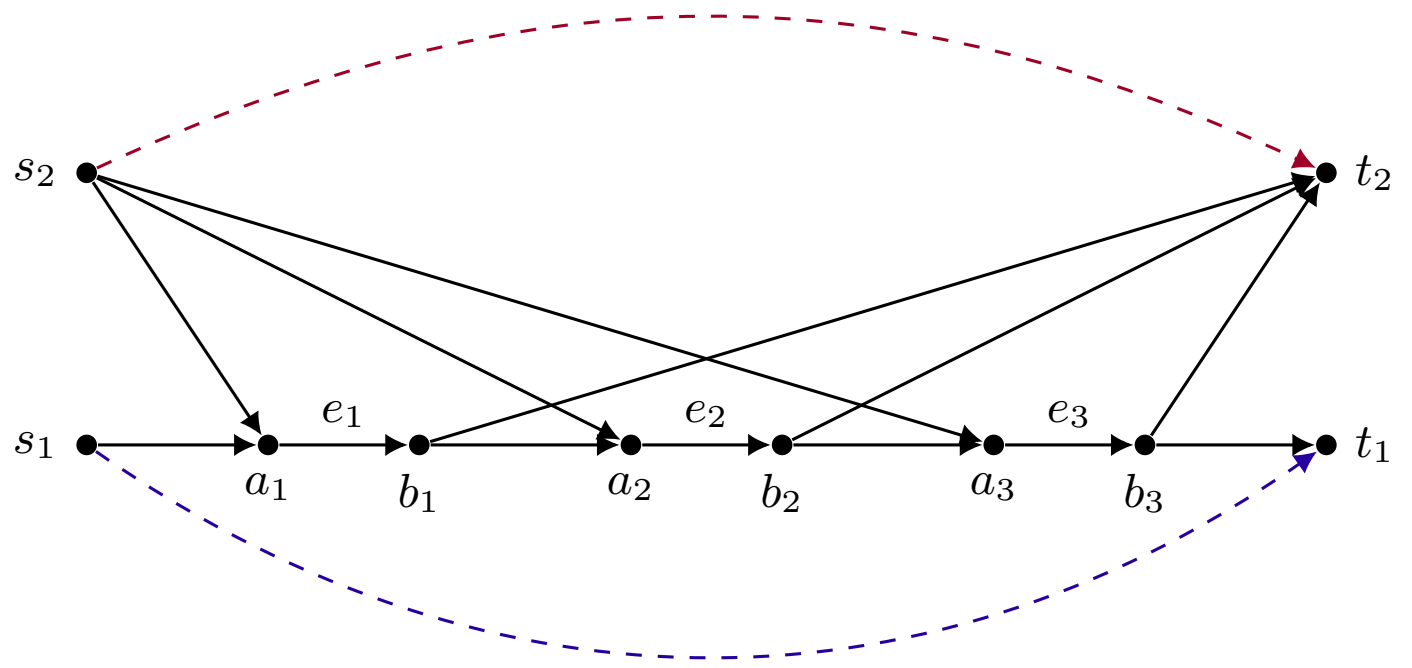

For $r\geq2$, define a topology $A_r$ on distinct vertices $s_1,s_2,t_1,t_2,a_1,b_1,\ldots,a_r,b_r$ and supply edges
\begin{align*}
  e_j=&(a_j,b_j),\ (s_2,a_j),\ (b_j,t_2)\ \text{for all } j\in[r], \\ 
  &(b_j,a_{j+1}) \text{ for all } j\in[r-1], \quad \text{ and } \quad
 (s_1,a_1),\ (b_r,t_1),
\end{align*}
see~\cref{gap:fig:Ar}. We call edges~$e_j$ the \emph{low edges} and all
other edges \emph{connectors}. The commodities are $(s_1,t_1)$ and
$(s_2,t_2)$. Commodity~$1$ has the unique path
$P=s_1,a_1,b_1,\ldots,a_r,b_r,t_1$, and commodity~$2$ can use any single
low edge or several consecutive ones. The connectors $(b_j,a_{j+1})$ are
needed: identifying~$b_j$ with~$a_{j+1}$ would give commodity~$2$ the
path $s_2,b_j,t_2$, which uses no low edge.

\begin{theorem}[Exact gaps for the forced-chain family]
\label{gap:thm:gapA}
For every $r\ge2$, $\beta_2(A_r)=\beta_3(A_r)=2-1/r$ and
$\beta_4(A_r)=1$.
\end{theorem}
\begin{proof}
\emph{Lower bound.} Give each low edge capacity~$1$ and each connector
capacity~$r$, and set $d:=(1,r-1)$. An edge set containing a connector
has capacity at least $r=d_1+d_2$. A nonempty proper subset of the low
edges blocks only commodity~$1$, and all low edges together block both
commodities. So~\eqref{C3}, and hence~\eqref{C2}, holds. In a routing within capacities~$\gamma u$, commodity~$1$ puts one unit
on each of the~$r$ low edges, and each of the $r-1$ units of
commodity~$2$ uses at least one low edge. The low edges have total
capacity~$\gamma r$, so $\gamma r\ge r+(r-1)$, that is,
$\gamma\ge2-1/r$. Conversely, $\gamma=2-1/r$ suffices: route
commodity~$1$ along~$P$, and split commodity~$2$ evenly over the~$r$
paths $s_2,a_j,b_j,t_2$, each carrying $(r-1)/r$. Each low edge then
carries $1+(r-1)/r=2-1/r$, and each connector carries at most one unit.

\emph{Upper bound.} Let weights satisfy~\eqref{C2}, and set
$\beta:=2-1/r$. As in the proof of~\cref{hier:unique}, every edge~$e$
of~$P$ forms a cut $\{e\}$ that blocks commodity~$1$, so $u_e\ge d_1$.
Hence commodity~$1$ can be routed on~$P$ within~$\beta u$. Every
commodity-$2$ path lies in $N:=G-\{s_1,t_1\}$. Let $F=\delta^+_N(X)$ be
an $s_2$--$t_2$ cut of~$N$, and let $q:=|F\cap P|$. Adding~$t_1$ to~$X$
shows that~$F$ is also a cut of~$G$. If $q=0$, then~\eqref{C2} gives
$u(F)\ge d_2$. If $q\ge1$, then~$F$ blocks both commodities, so
$u(F)\ge d_1+d_2$; moreover $u(F)\ge qd_1$ from the edges of~$P$ in~$F$.
The part of~$P$ in~$N$ has $2r-1$ edges, so it leaves~$X$ at most~$r$
times, and $q\le r$. Hence
\[
 qd_1+d_2=\frac{q-1}{q}\,(qd_1)+(d_1+d_2)
 \le\Bigl(2-\frac1q\Bigr)u(F)\le\beta u(F).
\]
So after commodity~$1$ is routed, every $s_2$--$t_2$ cut of~$N$ keeps
capacity at least~$d_2$ under~$\beta u$, and the max-flow/min-cut theorem
routes commodity~$2$. This proves the upper bound for~\eqref{C2}, and
hence for~\eqref{C3}. Finally, $\beta_4(A_r)=1$
by~\cref{hier:all-but-one}, since commodity~$1$ has a unique path.
\qed
\end{proof}

\subsection{The bypass family}
\label{gap:sec:B}

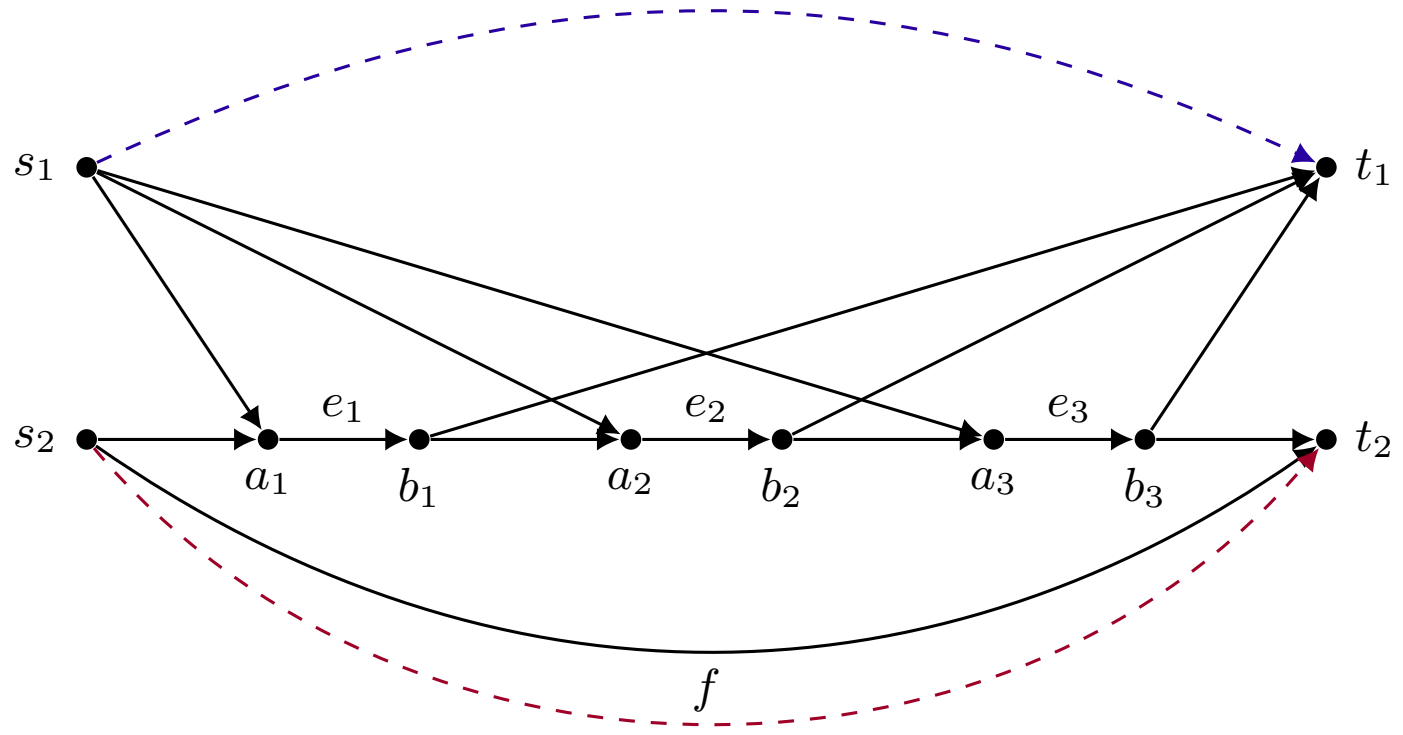
\begin{figure}[!b]
\centering
\begin{tikzpicture}[x=1cm,y=1cm]
  \foreach \i in {1,2,3}{
    \node[solidnode,label=below:$a_{\i}$] (a\i) at ({2.4*(\i-1)},0) {};
    \node[solidnode,label=below:$b_{\i}$] (b\i) at ({2.4*(\i-1)+1},0) {};
    \draw (a\i) edge[dedge] node[above]{$e_{\i}$} (b\i);
  }
  \node[solidnode,label=left:$s_2$]  (s2) at (-1.2,0) {};
  \node[solidnode,label=right:$t_2$] (t2) at (7,0) {};
  \node[solidnode,label=left:$s_1$]  (s1) at (-1.2,1.8) {};
  \node[solidnode,label=right:$t_1$] (t1) at (7,1.8) {};
  \draw (s2) edge[dedge] (a1) (b1) edge[dedge] (a2)
        (b2) edge[dedge] (a3) (b3) edge[dedge] (t2);
  \foreach \i in {1,2,3}{
    \draw (s1) edge[dedge] (a\i) (b\i) edge[dedge] (t1);
  }
  \draw (s2) edge[dedge,bend right=35] node[below]{$f$} (t2);
  \draw[demandedge,ccolor2] (s2) to[bend right=50] (t2);
  \draw[demandedge,ccolor1] (s1) to[bend left=25] (t1);
\end{tikzpicture}
\caption{The topology~$B_3$. Solid edges are supply edges; dashed edges
are demands. Crossings without a vertex are not junctions. In the
lower-bound weighting, the low edges have capacity~$1$, the edge~$f$
capacity~$0$, all other edges a large capacity~$U$, and
$(d_1,d_2)=(2,1)$.}
\label{gap:fig:Br}
\end{figure}

The topology~$B_r$ has the same vertices, the low edges
$e_j=(a_j,b_j)$, the connectors
\[
 (b_j,a_{j+1})\ (j\in[r-1]),\qquad (s_2,a_1),\ (b_r,t_2),\qquad
 (s_1,a_j),\ (b_j,t_1)\ (j\in[r]),
\]
and the edge $f=(s_2,t_2)$; see~\cref{gap:fig:Br}. Commodity~$1$ can use
any single low edge or several consecutive ones. Commodity~$2$ can
use~$f$ or the path $P=s_2,a_1,b_1,\ldots,a_r,b_r,t_2$. Both~$A_r$
and~$B_r$ are acyclic. They have $2r+4$ vertices, and $4r+1$ and $4r+2$
edges, respectively.

\begin{theorem}[Exact gaps for the bypass family]\label{gap:thm:gapB}
For every $r\ge2$, it holds that $\beta_3(B_r)=\beta_4(B_r)=2-1/r$.
\end{theorem}
\begin{proof}
\emph{Lower bound.} Set $d:=(r-1,1)$. Give each low edge capacity~$1$,
the edge~$f$ capacity~$0$, and each connector capacity
$U>(d_1+d_2)(|V|-1)$. We check~\eqref{C4}. If a binary length
vector~$\ell$ is~$1$ on some connector, its left-hand side is at
least~$U$, while every distance is at most $|V|-1$; so~\eqref{C4} holds.
Otherwise, let~$\eta$ be the number of low edges of length~$1$, and let
$\ell_f\in\{0,1\}$ be the length of~$f$. The left-hand side of~\eqref{C4} is~$\eta$. Commodity~$1$ can use any single low edge, so its distance
is~$1$ if all~$r$ low edges have length~$1$, and~$0$ otherwise.
Commodity~$2$ uses~$f$ or~$P$, so its distance is $\min\{\ell_f,\eta\}$. If
$\eta<r$, the right-hand side is $\min\{\ell_f,\eta\}\le\eta$. If $\eta=r$, it is
$(r-1)+\min\{\ell_f,r\}\le r=\eta$. So~\eqref{C4} holds.

In a routing within~$\gamma u$, commodity~$2$ cannot use~$f$, so it puts
one unit on every low edge, and each of the $r-1$ units of
commodity~$1$ uses at least one low edge. Hence $\gamma r\ge r+(r-1)$,
that is, $\gamma\ge2-1/r$. Conversely, $\gamma=2-1/r$ suffices: route
commodity~$2$ along~$P$, and split commodity~$1$ evenly over the~$r$
paths $s_1,a_j,b_j,t_1$, so that each low edge carries $1+(r-1)/r$.

\emph{Upper bound.} Let weights satisfy~\eqref{C3}, and set $c:=u_f$ and
$\beta:=2-1/r$. For every edge~$e$ of~$P$, the set $\{e,f\}$ blocks
commodity~$2$, so
\begin{equation}\label{gap:eq:bypass}
 u_e+c\ge d_2.
\end{equation}
We first route commodity~$2$ within~$\beta u$. If $\beta c\ge d_2$, send
all of it on~$f$. Otherwise, send $\beta c$ units on~$f$ and the
remaining $d_2-\beta c$ units along~$P$. This fits,
since~\eqref{gap:eq:bypass} gives
$\beta u_e\ge\beta(d_2-c)\ge d_2-\beta c$ for every edge~$e$ of~$P$.

Every commodity-$1$ path lies in $N:=G-\{s_2,t_2\}$. Let~$F$ be an
$s_1$--$t_1$ cut of~$N$. It blocks commodity~$1$, so $u(F)\ge d_1$
by~\eqref{C3}. We show that~$F$ keeps capacity at least~$d_1$ after
commodity~$2$ is routed; the max-flow/min-cut theorem then routes
commodity~$1$. If commodity~$2$ uses no edge of~$F$, this follows from
$\beta u(F)\ge u(F)\ge d_1$. Otherwise $\beta c<d_2$, and~$F$ contains
$q\ge1$ edges of~$P$, each carrying $d_2-\beta c$; moreover $q\le r$,
since the part of~$P$ in~$N$ has $2r-1$ edges and so leaves the vertex
set that defines~$F$ at most~$r$ times. Then $F\cup\{f\}$ blocks both
commodities, so $u(F)+c\ge d_1+d_2$ by~\eqref{C3}, and
summing~\eqref{gap:eq:bypass} over the~$q$ edges of~$P$ in~$F$ gives
$u(F)+qc\ge qd_2$. Therefore
\begin{align*}
 d_1+qd_2
 &\le\frac{q-1}{q}\bigl(u(F)+qc\bigr)+u(F)+c\\
 &=\Bigl(2-\frac1q\Bigr)u(F)+qc
 \le\beta\bigl(u(F)+qc\bigr),
\end{align*}
which says $\beta u(F)-q(d_2-\beta c)\ge d_1$. This proves the upper
bound under~\eqref{C3}, and hence under~\eqref{C4}.
\qed
\end{proof}

\subsection{Consequences}\label{gap:sec:consequences}

\begin{corollary}[No finite set of gap values]\label{gap:cor:allgaps}
For two commodities and each $\lambda\in\{1,2,3,4\}$, the supremum
of~$\beta_\lambda$ over all topologies is~$2$. For each $\lambda\in\{2,3,4\}$, the
values of~$\beta_\lambda$ include~$1$ and every $2-1/r$ with $r\ge2$, so no
level $\lambda \ge2$ has a finite set of gap values.
\end{corollary}
\begin{proof}
The upper bound is~\eqref{gap:upper}. At levels~$2$ and~$3$ the
family~$A_r$, and at levels~$3$ and~$4$ the family~$B_r$, give the
values $2-1/r$, which tend to two; level~$1$ is at least level~$2$.
Topologies with unique commodity paths have $\beta_\lambda=1$ for $\lambda \ge2$
by~\cref{hier:unique}.
\qed
\end{proof}

\begin{remark}[Integral gaps differ]\label{gap:rem:integer}
The integral gaps can be larger than the fractional ones. Take the
lower-bound weights on~$A_r$, which are integral, and capacities
$\gamma u$ with $\gamma<2$. In an integral routing, each low edge then
carries at most one unit. Commodity~$1$ already puts one unit on each of
them, so commodity~$2$ cannot be routed. Hence
$\beta_2^{\mathrm{int}}(A_r)=\beta_3^{\mathrm{int}}(A_r)=2$, whereas
$\beta_4^{\mathrm{int}}(A_r)=1$ by~\cref{hier:all-but-one,int:main}. In
the same way, the lower-bound weights on~$B_r$, where commodity~$2$ must
use every low edge, give $\beta_\lambda^{\mathrm{int}}(B_r)=2$ for $\lambda=2,3,4$;
these weights satisfy~\eqref{C4}, and hence also~\eqref{C2}
and~\eqref{C3}. So the supremum of~$\beta_\lambda^{\mathrm{int}}$ over
two-commodity topologies is~$2$ for every $\lambda\le4$.
\end{remark}

\section{Length tests with bounded coefficients}
\label{trunc:appendix}

For an integer $M\ge1$, let~\textup{(C5$_M$)} be~\eqref{C5} restricted
to length vectors contained in $\{0,1,\ldots,M\}^E$. Thus~\textup{(C5$_1$)}
is~\eqref{C4}, and~\textup{(C5$_M$)} becomes stronger as~$M$ grows. As
in~\cref{char:sec-cross}, a \emph{witness} is an integral length vector
$\ell\ge0$ that violates~\eqref{C5}; we write
$u\cdot\ell:=\sum_eu_e\ell_e$, so a witness satisfies
$u\cdot\ell<\sum_id_i\dist_\ell(s_i,t_i)$. An instance
satisfies~\textup{(C5$_M$)} exactly when it has no witness with all
lengths at most~$M$.

These are the truncated severance-value conditions of Onaga and
Kakusho~\cite{Onaga_Kakusho_71}. They showed that if~\textup{(C5$_M$)} holds, then all
demands can be routed after each is reduced by the same amount
$\sum_id_i\theta_i/M$, where the numbers~$\theta_i$ depend only on the
topology. They asked for the \emph{truncation number} of a given digraph,
the least~$M$ for which no reduction is needed. For a
topology~$(G,H)$, we let~$J(G,H)$ be the least~$M$ such that~$(G,H)$ is
\textup{(C5$_M$)}-sufficient, that is, every weighting
satisfying~\textup{(C5$_M$)} is feasible. Requiring positive capacities
does not change~$J(G,H)$: increasing capacities
preserves~\textup{(C5$_M$)}, and a limit of feasible instances is
feasible. Thus $J(G,H)=1$ if and only if $(G,H)\in\mathcal S_4$. We show
that~$J(G,H)$ is finite and at most exponential in the size of the
topology, but not bounded by any function of the number of commodities.

\begin{theorem}[Truncation numbers are finite]\label{gap:thm:hadamard}
Let $n:=|V|\ge2$ and $m:=|E|$. Every infeasible weighting has a
witness with all lengths at most $n^{(m+k-1)/2}$. Hence
$J(G,H)\le n^{(m+k-1)/2}$.
\end{theorem}
\begin{proof}
If a positive-demand pair is unreachable, then $\ell=0$ is a witness.
Otherwise, omit the zero-demand commodities; this does not
affect~\eqref{C5} and can only lower~$k$. Consider the cone
\[
 \mathcal K:=\bigl\{(\ell,q)\in\R_{\ge0}^{E}\times\R_{\ge0}^{k}:\
 q_i\le\ell(P)\text{ for all }i\in[k]\text{ and }P\in\mathcal P_i\bigr\},
\]
that is, the pairs $(\ell,q)$ with $q\ge0$ for which
$\sum_eu_e\ell_e\ge\sum_id_iq_i$ is valid for every feasible instance
(\cref{hier:appendix-recognition}). Write
$d\cdot q:=\sum_id_iq_i$. If $(\ell,q)\in\mathcal K$ and
$u\cdot\ell<d\cdot q$, then~$\ell$ is a witness, since
$\sum_id_i\dist_\ell(s_i,t_i)\ge d\cdot q$. Such a pair exists: take
any witness~$\ell$ and $q_i:=\dist_\ell(s_i,t_i)$.

The cone~$\mathcal K$ is defined by finitely many inequalities and lies
in the nonnegative orthant, so it is generated by its extreme rays. As
$u\cdot\ell<d\cdot q$ holds at some point of~$\mathcal K$, it holds on
some extreme ray. That ray is the set of points of~$\mathcal K$ at which
some $m+k-1$ linearly independent defining inequalities hold with
equality. The coefficient vectors of these inequalities have entries
in $\{0,1,-1\}$: an inequality $\ell_e\ge0$ or $q_i\ge0$ has one nonzero
entry, and an inequality $q_i\le\ell(P)$ has $|E(P)|+1\le n$ nonzero
entries. So each coefficient vector has Euclidean length at
most~$\sqrt n$. By Cramer's rule, the ray is spanned by the vector of
signed maximal minors of the matrix whose rows are these coefficient
vectors. By Hadamard's inequality, each such minor is at most
$(\sqrt n)^{m+k-1}$ in absolute value. Choosing the sign so that this
integral vector~$(\hat\ell,\hat q)$ lies on the ray gives
$u\cdot\hat\ell<d\cdot\hat q$, so~$\hat\ell$ is a witness with all
lengths at most $n^{(m+k-1)/2}$.

Finally, if $M\ge n^{(m+k-1)/2}$ and~\textup{(C5$_M$)} holds, there is
no witness with all lengths at most~$M$, so the instance is feasible.
\qed
\end{proof}

The bound depends only on the topology, not on~$u$ and~$d$, and the
number of paths does not enter it.

\begin{theorem}[Truncation numbers are unbounded]
\label{gap:thm:unboundedwitness}
For every $r\ge2$, the topology~$B_r$ of~\cref{gap:sec:B} has an
infeasible weighting with demands $(1,2)$ and positive integral
capacities that satisfies~\textup{(C5$_{r-1}$)}. Hence $J(B_r)\ge r$,
and for every $M\ge1$, \textup{(C5$_M$)} is strictly weaker
than~\textup{(C5$_{M+1}$)}, already on acyclic two-commodity topologies
with positive capacities.
\end{theorem}
\begin{proof}
Give each low edge and the edge~$f$ capacity~$1$, and each connector
capacity $U:=3r(2r+3)+1$. Commodity~$2$ can send only one unit on~$f$,
so it sends one unit along its path through all low edges, which fills
every low edge. Then commodity~$1$ cannot be routed, so the instance is
infeasible.

We check~\textup{(C5$_{r-1}$)}. A length vector with all lengths at
most $r-1$ that is positive on a connector has $u\cdot\ell\ge U$. Each
distance is at most $(r-1)(|V|-1)$, so the right-hand side is at most
$3(r-1)(2r+3)<U$. So consider length vectors that are zero on the
connectors. Let~$a$ be the smallest length of a low edge, and~$S$ the
total length of the low edges. The distance of commodity~$1$ is~$a$, and
that of commodity~$2$ is $\min\{S,\ell_f\}$. So the required inequality
is
\[
 S+\ell_f\ge a+2\min\{S,\ell_f\}.
\]
If $a=0$, this is $S+\ell_f\ge2\min\{S,\ell_f\}$. If $a\ge1$, then
$S\ge ra\ge r>\ell_f$, so the minimum is~$\ell_f$, and
$S+\ell_f-a-2\ell_f\ge(r-1)a-\ell_f\ge0$.

Length~$1$ on every low edge, length~$r$ on~$f$, and length~$0$ on
every connector give $u\cdot\ell=2r$ and right-hand side $1+2r$, so the
instance violates~\textup{(C5$_r$)}. Choosing $r=M+1$ gives the strict
separation.
\qed
\end{proof}

Bounded length tests also do not improve the worst-case gap. Define the
fractional gap of~\textup{(C5$_M$)} on a topology as
in~\cref{gap:appendix}, with~\textup{(C5$_M$)} in place
of~\textup{(C$\lambda$)}.

\begin{corollary}\label{gap:cor:truncgap}
For every fixed $M\ge1$, the supremum over two-commodity topologies of
the fractional gap of~\textup{(C5$_M$)} is two, even for acyclic
topologies with integral weights.
\end{corollary}
\begin{proof}
Take~$B_r$ with $r>M$. Give the low edges capacity~$1$, the edge~$f$
capacity~$0$, and the connectors capacity larger than
$M(|V|-1)(d_1+d_2)$; moreover, set $d:=(r-M,1)$. As in the previous proof,
length vectors that are positive on a connector satisfy the inequality,
and for the others it reads
\[
 S\ge(r-M)a+\min\{S,\ell_f\}.
\]
If $a=0$, it holds because $S\ge\min\{S,\ell_f\}$. If $a\ge1$, then
$S\ge ra>M\ge\ell_f$, so the minimum is~$\ell_f$, and
$S-(r-M)a-\ell_f\ge Ma-\ell_f\ge0$. Thus~\textup{(C5$_M$)} holds.

In a routing within~$\gamma u$, commodity~$2$ cannot use~$f$, so it puts
one unit on every low edge, and each of the $r-M$ units of
commodity~$1$ uses at least one low edge. Hence $\gamma r\ge r+(r-M)$.
Splitting commodity~$1$ evenly over the paths $s_1,a_j,b_j,t_1$ attains
$\gamma=2-M/r$, which tends to two as~$r$ grows. Since~\textup{(C5$_M$)}
implies~\eqref{C1}, routing the two commodities separately shows that
the gap is at most two.
\qed
\end{proof}

In summary, the worst-case truncation number grows at least linearly 
(\cref{gap:thm:unboundedwitness}) and at most exponentially
(\cref{gap:thm:hadamard}) in the size of the topology, and on each fixed
topology, \textup{(C5$_M$)} with $M\ge J(G,H)$ has gap one. Whether
truncation numbers can be super-polynomial in the size of the topology
remains open.

\section{Missing proofs from \cref{int:section}}
\label{apx:missing-proofs}

\begin{proof}[of \cref{int:unit}]
Each inequality in~\textup{(C$r$)} is $u(F)\ge b$ for some edge
set~$F$ and integer~$b$. In~\eqref{C1}-\eqref{C2}, $F$ is a
cut; in~\eqref{C3}, any edge set; in~\eqref{C4}, $F=\supp(
\ell)$, so $\sum_e u_e\ell_e=u(F)$.
The same $F$ also gives $x(F)\ge b$: every path of a separated or blocked commodity meets~$F$, and in~\eqref{C4} every $s_i$--$t_i$ path meets~$F$ at
least $\dist_\ell(s_i,t_i)$ times. Moreover, $b$ is finite since the existence of~$x$ makes all positive-demand pairs reachable, so no distance in the right-hand side
of~\eqref{C4} is infinite. If $e\in F$, then
$u(F)-x(F)\ge u_e-x_e>0$, hence $u(F)>b$ and, by integrality, $u(F)-1\ge b$. Inequalities with $e\notin F$ are unchanged.
\qed
\end{proof}

\begin{proof}[of \cref{int:total}]
Start with $c:=u$ and process the edges one at a time, keeping~$c$
integral and $(G,H,c,d)$ feasible. For the current edge~$e$, compute
\[
 \lambda_e:=\min\{x_e:x\text{ is a feasible multiflow for }(G,H,c,d)\}
\]
and a minimizer~$x$ by the edge-flow linear program
(\cref{sec:preliminaries}): shortcutting commodity cycles does not increase the objective $x_e$, so the path-flow and edge-flow linear programming optima coincide.

We claim that~$\lambda_e$ is integral. Otherwise, let~$c^\circ$
agree with~$c$ except that $c^\circ_e=\lceil\lambda_e\rceil$, which is at most~$c_e$.
Then~$x$ is feasible under~$c^\circ$ with slack
on~$e$, so by~\cref{int:unit}, the integral capacities
$c^\circ-\ind{\{e\}}$ satisfy~\eqref{C4} with demands~$d$.
The hypothesis therefore guarantees a feasible multiflow under these
capacities. It is also feasible under~$c$, but carries at most
$\lceil\lambda_e\rceil-1<\lambda_e$ on~$e$, a contradiction.
Thus,~$\lambda_e$ is integral. Set~$c_e:=\lambda_e$; the chosen
multiflow~$x$ remains feasible.

Let~$\tilde u$ be the final capacity vector; it is integral, at most~$u$, and admits a feasible multiflow.
Capacities only decrease, so a multiflow feasible under~$\tilde u$ was feasible when~$e$ was processed and $\lambda_e$ computed, so its flow on~$e$ is at least~$\lambda_e$.
Its flow is also at most $\tilde u_e=\lambda_e$, because the capacity
of~$e$ was not changed after its processing step. Hence every feasible
multiflow under~$\tilde u$ saturates every edge.
The procedure solves one linear program per edge.
\qed
\end{proof}

\end{document}